\documentclass[11pt]{amsart}

\newcounter{nootje}
\usepackage{amssymb,amsmath,enumerate,epsfig}
\usepackage{amsfonts,xypic}
\usepackage{a4,latexsym,parskip}
\usepackage{hyperref}
\hypersetup{pdftitle=Fermat}

\newcommand\rk{\mathop{\rm rk}}
\newcommand\tors{{\text{torsion}}}

\newcommand{\C} [1][]{\mathbb{C}^{#1}}
\newcommand{\Q} [1] []{\mathbb{Q}_{#1}}
\newcommand{\N} [1][] {\mathbb{N}_{#1}}

\newcommand{\F}{\mathbb{F}}
\newcommand{\Z}{\mathbb{Z}}
\renewcommand\P{\mathbb{P}}

\newcommand{\OO}{\mathcal{O}}

\newcommand{\e}{\varepsilon}

\newcommand{\PP}{\mathbb{P}}
\newcommand{\NS}{\mathop{\rm NS}}
\newcommand{\Num}{\mathop{\rm Num}}
\newcommand{\A}{\mathfrak{A}}
\newcommand{\B}{\mathfrak{B}}
\newcommand{\D}{\mathfrak{D}}
\newcommand\disc{\mathop{\rm disc}}
\newcommand\Hom{\mathop{\rm Hom}}
\newcommand{\li}{\mathfrak{l}}
\newcommand\ptors{{\text{$p$-torsion}}}

\newtheorem{Theorem}{Theorem}[section]
\newtheorem{Proposition}[Theorem]{Proposition}
\newtheorem{Lemma}[Theorem]{Lemma}
\newtheorem{Criterion}[Theorem]{Criterion}
\newtheorem{Result}[Theorem]{Result}
\newtheorem{Corollary}[Theorem]{Corollary}

\theoremstyle{remark}
\newtheorem{Remark}[Theorem]{Remark}

\theoremstyle{definition}

\begin{document}
\setlength{\unitlength}{1cm}

\title{Lines on Fermat surfaces}

\subjclass[2000]{Primary 14J25; 
Secondary 11G25, 14C22.}

\keywords{Fermat surface, N\'eron-Severi group, supersingular reduction}

\author{Matthias Sch\"utt}
\address{Institute for Algebraic Geometry, Leibniz University Hannover, Welfengarten 1, 30167 Hannover, Germany}
\email{schuett@math.uni-hannover.de}
\urladdr{http://www.iag.uni-hannover.de/\~{}schuett/}

\author{Tetsuji Shioda}
\address{Department of Mathematics, Rikkyo University, Tokyo 171-8501, Japan, and\;\;\;
RIMS,
Kyoto University,
Kyoto 606-8502,
Japan}
\email{shioda@rikkyo.ac.jp}
\urladdr{http://www.rkmath.rikkyo.ac.jp/math/shioda/}


\author{Ronald van Luijk}
\address{Mathematisch Instituut, Universiteit Leiden, Postbus 9512, 2300 RA, Leiden, The Netherlands}
\email{rvl@math.leidenuniv.nl}
\urladdr{http://www.math.leidenuniv.nl/\~{}rvl}

\thanks{Partial funding from DFG under grant Schu 2266/2-2 and JSPS under Grant-in-Aid for Scientific Research (C) No.~20540051 is gratefully acknowledged.}

\date{September 30, 2009}

\begin{abstract}
We prove that the N\'eron-Severi groups of several complex Fermat surfaces are generated by lines. 
Specifically, we obtain these new results for all degrees up to $100$ that are relatively prime to $6$. The proof uses reduction modulo a supersingular prime. The techniques are developed in detail. 
They can be applied to other surfaces and varieties as well.
\end{abstract}

\maketitle

\section{Introduction}

Fermat varieties have been a classical object of study in geometry and arithmetic. Here we consider the smooth projective surface of degree $m\in\N$
\[
S:\;\;\;\{x_0^m+x_1^m+x_2^m+x_3^m=0\}\subset\PP^3.
\]
This paper is concerned with the N\'eron-Severi group $\NS(S)$ of $S$ over the complex numbers, consisting of divisors up to algebraic equivalence. 

In general, it is hard to compute the N\'eron-Severi group of a variety. 
The cohomology of Fermat surfaces, however, admits a decomposition into eigenspaces with 
respect to an abelian subgroup of the automorphism group. Combinatorial data 
give the Picard number $\rho(S)$, the rank of $\NS(S)$. 
A rational basis of $\NS(S)$ (i.e.~a basis of $\NS(S)\otimes\Q$) was determined in \cite{AS}
up to certain cycles induced from Fermat surface with degree $m$ in the range $12\leq m\leq 180$.

The cycles exhibited in \cite{AS} involve some particularly prominent divisors on $S$, namely the $3m^2$ obvious lines. The lines generate $\NS(S)$ rationally if and only if $m\leq 4$ or $(m,6)=1$. 
In Proposition~\ref{Prop:basis}, we will improve the results from \cite{AS} in the sense that we identify a rational basis consisting of lines explicitly. 

A natural question now is in which of the above cases the lines generate the full N\'eron-Severi group. 
As opposed to \emph{rational generation}, we refer to this property as \emph{integral generation}.
Integral generation is known to hold true for $m\leq 4$, as we will review in section \ref{s:rat}. Our main result is the following:

\begin{Theorem}\label{thm}
Let $m\leq 100$ be a positive integer. Then the N\'eron-Severi group of the complex Fermat surface $S$ of degree $m$ is integrally generated by lines if and only if $m\leq 4$ or $(m,6)=1$. 
\end{Theorem}

We shall use supersingular reduction to prove the theorem.
The technique is briefly outlined below; a full account will be given in section~\ref{s:tech}.
For the degrees $m<17$, the method is applied separately in sections \ref{s:5}--\ref{s:13} to exhibit a proof of the corresponding part of Theorem~\ref{thm}. 
In section \ref{s:out}, we develop an extension of the supersingular reduction technique that is less  involved computationally.
This technique is applied to the remaining degrees in section \ref{ss:Fermat} to complete the proof of Theorem \ref{thm}.

We give a brief outline of the supersingular reduction technique.
Starting from a complex Fermat surface $S$, we consider the reduction $S_p$ modulo a good prime $p$.
By choosing a supersingular prime, we achieve good control of the discriminant of the N\'eron-Severi lattice of the reduction $S_p$ (Theorem \ref{Thm:NS}).
Then we compare the discriminants of two lattices:
on the one hand the  sublattice of finite index in $\NS(S)$ generated by lines, 
on the other hand a suitable (often finite-index) 
sublattice of $\NS(S_p)$ where we complement the reductions of the original lines by some divisors that are peculiar to the chosen characteristic (cf.~section \ref{s:add}).
Unless these discriminants have a common square factor, this method suffices to prove that the sublattice generated by lines is already the full N\'eron-Severi lattice by Criterion \ref{Crit}.

In spirit the supersingular reduction technique is related to a method to compute the Picard number of a projective surface which was introduced by one of the authors in \cite{vL}.
Namely it was proved that certain K3 surfaces have Picard number one by reducing modulo two different primes. 
From the Lefschetz fixed point formula, one would derive that the reductions had Picard number (at most) two. 
Then one would find divisors peculiar to the respective characteristic and compare the resulting discriminants of the N\'eron-Severi lattices. Once they did not match up to a square factor, it would follow that the original surface had Picard number one.

The supersingular reduction technique compares sublattices of $\NS(S)$ and $\NS(S_p)$ for a supersingular prime $p$, 
while the method in \cite{vL} required two suitable reductions.
Both methods are greatly inspired by the Tate conjecture \cite{TT},
and in fact the equivalent statement of the Artin-Tate conjecture \cite{Milne} plays a crucial role for several aspects (cf.~Theorem \ref{Thm:NS} and \cite{Klo}).

The computations were carried out with MAGMA.
Programs and scripts are available from the third author's webpage.
We are indepted to Bas Edixhoven for the use of his computer.

\section{Preliminaries on projective surfaces and lattices}\label{s:prelim}

In this section, we recall some basic facts about lattices, projective surfaces and divisors that are relevant for our purposes.
In view of Fermat surfaces, we will mostly be concerned with smooth surfaces in $\PP^3$.
For general background, the reader might confer \cite{BHPV} or \cite{Mumford}.

Throughout this paper, every lattice is assumed to be integral unless otherwise stated. In other words, a lattice is a finitely generated free abelian group $\Lambda$, together with a symmetric bilinear pairing $\langle  \cdot, \cdot \rangle  \colon \Lambda \times \Lambda \to \Z$ that is nondegenerate, i.e., the induced map $\Lambda \to \Hom(\Lambda,\Z)$ is injective. 
The discriminant of a lattice $\Lambda$ is the determinant of the Gram matrix 
$\big(\langle  x, y \rangle \big)_{x,y}$, where $x$ and $y$ run through any chosen basis of $\Lambda$; the discriminant is independent of the choice of basis. If $L$ is a finite-index sublattice of a lattice $\Lambda$, then their discriminants are related through the equality 
\[
\disc (L) = [\Lambda : L ]^2 \cdot \disc (\Lambda).
\]
We say that a sublattice $L$ is primitive in $\Lambda$ if the quotient $\Lambda/L$ is torsion-free.
This is only possible if $L$ has positive corank in $\Lambda$ or $L=\Lambda$.

If $k$ is a field and $L$ a lattice, then $L_k$ denotes the vector space $L\otimes_{\Z} k$.
For any vector space $V$ over a field $k$, we denote its dual $\Hom(V,k)$ by $V^*$.

On any projective surface $X$, the curves generate the group Div$(X)$ freely.
This group can be endowed with a meaningful structure by dividing out by some equivalence relation such as linear equivalence $\sim$, algebraic equivalence $\approx$ or numerical equivalence $\equiv$ (with implications from left to right).

Two curves are algebraically equivalent if they move within a family of divisors on $X$ over some curve (for instance any fibration has algebraically equivalent fibers).
The N\'eron-Severi group of a projective surface $X$ is defined as the quotient
\[
 \NS(X) = \mbox{Div}(X)/\approx.
\]
Its rank is called the Picard number, denoted by $\rho(X)$.
The N\'eron-Severi group depends on the chosen base field of the variety (such as number fields, finite fields).
In this paper, we are concerned with geometric invariants;
hence we always consider the geometric N\'eron-Severi groups,
i.e.~for a base change of the surface to an algebraic closure of its base field ($\C, \bar\Q, \bar\F_p$).
Whenever $X$ is a surface over $\C$ and we want to reduce it modulo a prime $p$, it is implicitly understood that 
we first take an integral model of $X$ that has good reduction at $p$; 
the surface in the reduction can then be considered over $\bar\F_{p}$.

Two divisors are numerically equivalent if they return the same intersection numbers with all divisors on $X$ 
(or equivalently with all divisor classes in $\NS(X)$).
The corresponding quotient is denoted by $\Num(X)$.
It is known 
that the only difference between algebraic and numerical equivalence 
lies in the torsion in $\NS(X)$:
\[
 \Num(X) = \NS(X)/\text{torsion}.
\]
In particular, these notions coincide if $X$ is (algebraically) simply connected.
This holds for large classes of varieties such as complete smooth intersection in $\PP^n$ of dimension greater than one.
In consequence, for any smooth surface $X$ in $\PP^3$, the N\'eron-Severi group is torsion-free.
The intersection form endows $\NS(X)$ with the structure of a lattice, also called the N\'eron-Severi lattice.
By the Hodge index theorem, the N\'eron-Severi lattice has signature $(1,\rho(X)-1)$.

We have seen that it suffices to compute intersection numbers to understand the N\'eron-Severi groups of Fermat surfaces.
Self-intersection numbers involve a subtlety as they can be negative, depending on the chosen surface.
For a (smooth) irreducible curve $C$ on a surface $X$, one can compute $C^2$ through the adjunction formula:
\[
 2g(C)-2 = C^2 + C.K_X.
\]
Here $g(C)$ is the genus of $C$ and $K_X$ denotes the canonical divisor of $X$.
Often the canonical divisor can be expressed through a hyperplane section $H$.
For a smooth surface of degree $m$ in $\PP^3$, one has $K_X=(m-4)H$.
For a line $\li$ and a conic $Q$ (both rational curves, thus of genus zero), one obtains the following self-intersection numbers on such a surface $X$:
\[
 \li^2 = 2-m,\;\;\; Q^2 = 6-2m.
\]
More generally, if $C$ is a rational curve of degree $d$ on a smooth surface of degree $m$ in $\PP^3$, then 
$C^2 = 4d-2-dm$.

We conclude this section by indicating how to compute the Betti numbers and Hodge numbers of a smooth (complex) surface $X$ of degree $m$ in $\PP^3$.
We have already mentioned that $b_1(X)=q(X)=0$.
By Serre duality, the geometric genus equals 
\[
p_g(X)=h^2(X,\OO_X)=h^0(X,K_X)=h^0(\PP^3, \OO_{\PP^3}(m-4))=\binom{m-1}3.
\]
Thus we compute the Euler characteristic $\chi(\OO_X)=h^0(\OO_X(X))-q+p_q(X) = 1+p_g(X)$.
The topological Euler number $e(X)$ (which can be defined as the alternating sum of Betti numbers in arbitrary characteristic) can be computed by Noether's formula
\[
 12\chi(\OO_X) = e(X) + K_X^2.
\]
Here $K_X^2 = m(m-4)^2$. 
Then the second Betti number is calculated as $b_2(X)=e(X)-2$, as we have $b_0=b_4=1$ and $b_1=b_3=0$ by Poincar\'e duality. One finds 
\[
b_2(X) = m^3-4m^2+6m-2.
\]
Over $\C$, we obtain the Hodge number $h^{1,1}=b_2(X)-2p_g(X)$.
The Picard number relates to these invariants as follows:
\begin{itemize}
 \item 
in characteristic zero, $\rho(X)\leq h^{1,1}(X)$ by Lefschetz' theorem;
\item
in positive characteristic, $\rho(X)\leq b_2(X)$ by Igusa's theorem.
\end{itemize}
Surfaces attaining an equality in the latter setting are often called \emph{supersingular}.
We will recall some of their properties in section \ref{s:tech} and use them for our supersingular reduction technique.

\section{Rational generation of $\NS$}
\label{s:rat}

The cohomology of Fermat varieties admits a decomposition into eigenspaces with respect to an abelian subgroup of the automorphism group. 
According to work by Katz and Ogus, it splits into one-dimensional eigenspaces; 
we review these concepts below starting with (\ref{eq:mu}). 
It is well known which eigenspaces are algebraic, and in the surface case, even which eigenspaces correspond to lines. 

\begin{Theorem}[Shioda {\cite{Sh-Fermat}}]
\label{Thm:AS}
Let $S$ denote the complex Fermat surface of degree $m$.
The $\Q$-vector space $\NS(S)\otimes_{\Z}{\Q}$ is generated by divisor classes of lines if and only if $m\leq 4$ or $m$ is coprime to $6$.
\end{Theorem}

Before reviewing the proof of the theorem,
we comment on the main problem of this paper whether, for the appropriate degrees, lines generate $\NS(S)$ fully or only up to finite index. We now review the current knowledge about this problem. 

For $m\leq 3$, the generation problem has a positive answer. These Fermat surfaces are rational. 
For $m=1, 2$, the statement is almost trivial, corresponding to $\PP^2$ and $\PP^1\times\PP^1$.
Any smooth projective cubic complex surface contains 27 lines.
Their configuration has been  studied in great detail.
In fact, any smooth cubic surface is isomorphic to the projective plane $\PP^2$ blown up in six distinct  points.

For $m=4$, the K3 case, the answer was conjectured to be positive, but unknown until Mizukami in 1975 proved the affirmative \cite{Mizukami}. We will review the history of the original proof and provide an alternative proof using our technique of supersingular reduction in section \ref{s:4}.
Our Theorem~\ref{thm} provides the first answer to the question for Fermat surfaces of general type.

In the sequel we shall sketch the line of argument from \cite{Sh-Fermat} 
for later use in the next section.
In order to prove Theorem \ref{Thm:AS}
it clearly suffices to prove the corresponding statement for $\NS(S)\otimes\C$.
Hence we will mostly work with the latter vector space in this section and analyse when it is generated by lines.

First we fix notation for the $3\,m^2$ lines on $S$, the Fermat surface of degree $m$.
Throughout the paper, we denote by $\mu_n$ the group of $n$-th roots of unity over a given field.
Let $\omega\in\mu_{2m}$ such that $\omega^m=-1$. Then for any $\zeta,\eta\in\mu_m$ we have the lines 
\begin{eqnarray*}
\mathfrak l_1(\zeta, \eta) & = & \{[\lambda, \omega\,\zeta \lambda, \mu, \omega\,\eta\,\mu];\;\; [\lambda,\mu]\in\PP^1\},\\
\li_2(\zeta, \eta) & = & \{[\lambda,  \mu, \omega\,\zeta \lambda, \omega\,\eta\,\mu];\;\; [\lambda,\mu]\in\PP^1\},\\
\li_3(\zeta, \eta) & = & \{[\lambda, \mu, \omega\,\eta\,\mu, \omega\,\zeta \lambda];\;\; [\lambda,\mu]\in\PP^1\}.
\end{eqnarray*}
On $S$, the abelian group $\mu_m^4/\mu_m$ acts by  multiplication on homogeneous coordinates:
\begin{eqnarray}\label{eq:mu}
\;\;\;\;\;
g=[\zeta_1,\zeta_1, \zeta_2, \zeta_3]\in\mu_m^4/\mu_m:\;\;\; [x_0,x_1,x_2,x_3] \mapsto [\zeta_0\,x_0,\zeta_1\,x_1, \zeta_2\,x_2,\zeta_3\,x_3].
\end{eqnarray}
The character group of $\mu_m^4/\mu_m$ is isomorphic to the kernel of the map 
$$
\textstyle{\sum}\colon (\Z/m\Z)^4 \to \Z/mZ, \qquad \alpha = (a_0,a_1,a_2,a_3) \mapsto \sum_i a_i,
$$
where $\alpha$ sends $g=[\zeta_0,\zeta_1,\zeta_2,\zeta_3] \in \mu_m^4/\mu_m$ to $\alpha(g) = \prod_i \zeta_i^{a_i} \in \mu_m$. 
We shall consider the eigenspaces of $H^2(S)$ for the induced action of $\mu_m^4/\mu_m$ with character $\alpha$ in the 
following subset of the character group
\[
\A_{m}:=
\{\alpha=(a_0,a_1,a_2,a_3)\in \ker \Sigma \,\,|\,\, a_i \neq 0 \,\}.
\]
For $\alpha\in\A_{m}$, the corresponding eigenspace $V(\alpha)\subset H^2(S)$ with character $\alpha$ is defined by the condition
that $g^*|_{V(\alpha)}$ acts as multiplication by $\alpha(g)$ for all $g \in \mu_m^4/\mu_m$.
By results of Katz \cite[\S 6]{Katz} and Ogus \cite[\S 3]{Ogus} (which hold more generally true for Fermat varieties of any dimension), each $V(\alpha)$ is one-dimensional, and
\[
H^2(S) = V_0 \oplus \bigoplus_{\alpha\in\A_{m}} V(\alpha).
\]
Here $V_0$ corresponds to the trivial character and is spanned by the hyperplane section.
One easily checks that $\#\A_m = (m-1)(m^2-3m+3)$, so that indeed $\# \A_m+1 = b_2(S)$.

Up to this point, the whole argument does not depend on the characteristic and works for any appropriate cohomology theory.
From now on, we specialise to the complex case.
Writing $\alpha = (a_0,\hdots,a_3)\in\A_m$ with canonical representatives $0<\tilde{a}_i<m$, we define
\[
|\alpha| = (\tilde{a}_0+\hdots+\tilde{a}_3)/m.
\]
Then the eigenspace $V(\alpha)$ has the Hodge weights $(|\alpha|-1, 3-|\alpha|)$.
In order to decide whether $V(\alpha)$ is algebraic, 
we let $(\Z/m\,\Z)^*$ operate on $\A_{m}$ coordinatewise by multiplication. 
As a consequence of Lefschetz' theorem, $V(\alpha)$ is algebraic
if and only if every element in the $(\Z/m\,\Z)^*$-orbit of $\alpha$ has Hodge weight $(1,1)$,
i.e., if and only if $|r\alpha|=2$ for all $r \in (\Z/m\Z)^*$.

To collect the corresponding $\alpha$, we define the subset $\mathfrak{B}_{m}\subset\A_{m}$ as follows:
\[
\alpha\in\mathfrak{B}_m\;\; \Longleftrightarrow \;\; \forall\, r\in(\Z/m\,\Z)^*: \; |r\alpha|=2.
\]
The space $V(\alpha)$ is algebraic if and only if $\alpha\in\mathfrak{B}_m$. Hence
\[
\rho(S) = \#\mathfrak{B}_m+1.
\]
By \cite{Sh-Fermat}, the span of the lines  is also known: In $\NS(S)\otimes\C$, this is 
\begin{equation}\label{VD}
V_0 \oplus \bigoplus_{\alpha\in \D_m} V(\alpha),
\end{equation}
where $\D_m\subseteq\B_m$ denotes the subset of decomposable elements $\alpha$, i.e.,~those 
$\alpha \in \B_m$ for which there is some index $j>0$ such that $a_0+a_j=0$. 
Then one easily computes
\begin{eqnarray}\label{eq:D_m}
 \# \D_m = 3\,(m-1)\,(m-2) + 
\begin{cases}
 0, & \text{ if $m$ is odd},\\
1, & \text{ if $m$ is even}.
\end{cases}
\end{eqnarray}
We now recall why the lines generate the space in (\ref{VD}). 
This will be achieved by establishing a $\C$-linear combination of lines which is a non-zero eigendivisor for the character $\alpha\in\D_m$. 

More specifically, let $\D_m^j$ denote the subset of decomposable elements in $\D_m$ such that $a_0+a_j=0$. 
Note that $\D_m^j\cap \D_m^k\neq \emptyset$ for all $1\leq j,k \leq 3$ -- a fact 
that will be crucial to our later analysis of an explicit basis of lines.
Depending on $j$, we give an eigendivisor with character for each $\alpha\in\D_m^j$:
\begin{eqnarray*}
\alpha\in\D_m^1:\;\; &
w_1(\alpha) = \sum_{\zeta,\eta} \zeta^{a_1}\,\eta^{a_3} \li_1(\zeta,\eta),\\
\alpha\in\D_m^2:\;\; &
w_2(\alpha) = \sum_{\zeta,\eta} \zeta^{a_2}\,\eta^{a_3} \li_2(\zeta,\eta),\\
\alpha\in\D_m^3:\;\; &
w_3(\alpha) = \sum_{\zeta,\eta} \zeta^{a_3}\,\eta^{a_2} \li_3(\zeta,\eta),
\end{eqnarray*}
where the sum is over all $\zeta,\eta \in \mu_m$. 
By construction, almost all of these eigendivisors are orthogonal:
\begin{eqnarray}\label{eq:ortho}
w_i(\alpha).H=0, \;\;\; w_i(\alpha).w_j(\beta)=0 \;\;\; \text{ if } \alpha\neq-\beta \;\;\;\;\;\; (i,j=1,2,3).
\end{eqnarray}
which is easily computed thanks to the following intersection behaviour:
\begin{eqnarray}\label{eq:inter}
\li_i(\zeta,\eta).\li_j(\zeta',\eta') \neq 0 \Leftrightarrow 
\begin{cases}
\zeta=\zeta' \text{ or } \eta=\eta', & i=j,\\
\zeta\,\eta'=\zeta'\,\eta, & (i,j)=(1,2),\\
\zeta'=\omega^2\,\eta\,\zeta\,\eta', & (i,j)=(1,3),\\
\zeta\,\eta = \zeta'\,\eta', & (i,j)=(2,3).
\end{cases}
\end{eqnarray}
From
the intersection number
\begin{eqnarray}\label{eq:inter-w}
w_j(\alpha).w_j(-\alpha)=-m^3
\end{eqnarray}
it follows that $w_j(\alpha)\neq 0$. 
We conclude that $V(\alpha) \subset \NS(S)\otimes\C$ is contained in the span of the lines. Denote this span by $L$. 
Clearly, also $H$ and thus $V_0$ can be expressed by lines (cf.~(\ref{rel:H}), (\ref{rel:H'})), so we derive 
the inclusion $\subset$ of the following equality
\begin{equation}\label{VDinL}
V_0 \oplus \bigoplus_{\alpha\in \D_m} V(\alpha) = L.
\end{equation}
The other inclusion follows from the fact that every line can be expressed in terms of $H$ and the $w_j(\alpha)$ for $\alpha\in\D_m$ (cf.~\cite[(17)]{Sh-Fermat}).
In particular, we have
\[
\mbox{rank}(L) = 1+ \# \D_m.
\]
\begin{proof}[Proof of Theorem~\ref{Thm:AS}]
We have seen that the span of lines $L$ has rank $1+\#\D_m$. 
On the other hand, $\rho(S)=1+\#\B_m$.
From \cite[Theorem~6]{Sh-Fermat} we know that 
\[
\D_m=\B_m\;\;\; \Longleftrightarrow \;\;\; m\leq 4\;\;\text{ or }\;\; (m,6)=1.
\]
This proves that the lines generate $\NS(S)\otimes \C$ exactly in the cases of Theorem~\ref{Thm:AS}. The corresponding statement for 
$\NS(S)\otimes \Q$ follows. 
\end{proof}

\begin{Corollary}\label{Cor:discr-m}
The lattice $\Lambda$ generated by the lines has discriminant dividing $m^r$ for $r=3\#\D_m+1$.
\end{Corollary}

\begin{proof}
Consider the $\Z[\zeta_m]$-lattice $\Lambda\otimes\Z[\zeta_m]$.
It contains the  finite-index sublattice $\Lambda'$ generated by $H$ and the $w_j(\alpha)$ for $j =1,2,3$ and $\alpha\in\D_m^j$.
The given generators of $\Lambda'$ have intersection matrix $Q'$ of determinant $m^r$ for $r=3\,\# \mathfrak{D}_m+1$ by (\ref{eq:ortho}) and (\ref{eq:inter-w}).
The determinant of $Q'$ equals the discriminant of $\Lambda$ times a square in $\Z[\zeta_m]$ (the square of the determinant of the matrix in $M_\rho(\Z[\zeta_m])\cap GL_\rho({\Q}[\zeta_m])$ that expresses the given basis of $\Lambda'$ in terms of a basis of $\Lambda$).
Hence $\Lambda$ has discriminant that divides $m^r$.
\end{proof}

\section{Rational basis of lines}
\label{s:lines-basis}

In this section, we will work out an explicit rational basis of the lattice $L$ generated by the lines in $\NS(S)$ for the complex Fermat surface $S$ of degree $m$. 
For this, we fix another notation for the lines.
Since we are concerned with odd $m$, we can set $\omega=-1$. 
Then we fix a primitive $m$-th root of unity $\gamma$. We introduce the short-hand notation
\[
\li_j(\gamma^k, \gamma^l) = \li_j(k,l)
\]

\begin{Proposition}[Rational basis for $m$ coprime to $6$]
\label{Prop:basis}
Assume that $(m,6)=1$ and that the ground field has characteristic zero.
Then the following lines form a basis of $\NS(S)\otimes \Q$:
\[
\mathcal{B} = \{\li_j(k,l);\;\; j=1,2,3,\; 0\leq k<m-1,\; 0< l<m-1\} \cup \{\li_1(m-1,1)\}
\]
\end{Proposition}

\emph{Proof:}
We shall use relations between lines and the hyperplane class $H$. Clearly 
\begin{eqnarray}
\label{rel:H}
H & = & \sum_\zeta \li_i(\zeta,\eta)\\
\label{rel:H'}
& = & \sum_\eta \li_j(\zeta,\eta)
\end{eqnarray}
for any fixed $\eta$ resp.~$\zeta$ and independent of the index. Taking the sum of the lines $\li_1(\cdot,1)$, we see that $H$ is in the span of $\mathcal{B}$. In consequence, all $\li_i(m-1,l)$ for $1<l<m-1$ can be expressed by $\mathcal{B}$ as well. It remains to write the lines $\li_i(\cdot,0), \li_j(\cdot,m-1)$ in terms of the previous lines.

A second set  of relations is derived for all those $\alpha\in\D_m^i\cap\D_m^j$ for some $i\neq j$. Since $V(\alpha)$ is always one-dimensional, we have
\[
V(\alpha) = \C\,w_i(\alpha) = \C\,w_j(\alpha),
\]
so the two eigendivisors are multiples of each other. Recall that each eigendivisor $w_j(\alpha)$ intersects its complex conjugate $w_j(-\alpha)$ with intersection multiplicity $-m^3$.

\textbf{Claim:} 
Let $i\neq j$ and $\alpha\in\D_m^i\cap\D_m^j$. Then
\begin{eqnarray}\label{claim}
w_i(\alpha) = -w_j(\alpha).
\end{eqnarray}

Recall the orthogonality for eigendivisors with character from (\ref{eq:ortho}). To see the claim, it thus suffices to compute the intersection number
\[
w_i(\alpha).w_j(-\alpha) = m^3.
\]
This is easily verified thanks to the intersection behaviour of the lines in (\ref{eq:inter}).

The coefficients of the lines in the relations (\ref{claim}) involve $m$-th roots of unity.
In order to derive relations over $\Q$, we shall now simplify the above relations by multiplying with fixed powers of a varying root $\varepsilon\in\mu_m$.

For any pair $(i,j)$ with $i\neq j$,
we define the map
\begin{eqnarray*}
\alpha_{i,j}:\Z/m\Z - \{0\} & \to & \D_m^i\cap\D_m^j\\
r \;\;\;\;\;\;\; & \mapsto & \;\; \alpha_{i,j}(r)
\end{eqnarray*}
by setting $a_0=r$. Then $a_i=a_j=-r$ and $a_k=r$ with $\{i,j,k\}=\{1,2,3\}$.
For any $\e\in\mu_m$ and $(i,j)$ with $i\neq j$, we then consider the relations of divisors obtained from (\ref{claim})
\[
\sum_{r\in\Z/mZ-\{0\}} \varepsilon^r\,w_i(\alpha_{i,j}(r)) = - \sum_{r\in\Z/mZ-\{0\}} \varepsilon^r\,w_j(\alpha_{i,j}(r)).
\]
Both sums simplify greatly. For instance, 
\begin{eqnarray*}
\sum_{r\in\Z/mZ-\{0\}} \varepsilon^r\,w_1(\alpha_{1,2}(r)) & = & \sum_{\zeta, \eta}\,\sum_r \left(\frac{\e\,\eta}\zeta\right)^r\,\li_1(\zeta,\eta)\\
& = & (m-1)\,\sum_{\zeta=\e\,\eta} \li_1(\zeta, \eta) - \sum_{\zeta\neq\e\,\eta} \li_1(\zeta, \eta)\\
& \stackrel{(\ref{rel:H})}{=} & m\,\left( \sum_{\zeta=\e\,\eta} \li_1(\zeta, \eta) - H\right).
\end{eqnarray*}
Analogous sums for the other indices result in the following $3\,m$ relations (depending on the choice of $\e\in\mu_m$):
\begin{eqnarray}
\label{rel:12}
\sum_{\zeta=\e\,\eta} \li_1(\zeta, \eta)  & = & - \sum_{\zeta=\e\,\eta} \li_2(\zeta, \eta)\\
\label{rel:13}
\sum_{\zeta\,\eta=\e} \li_1(\zeta, \eta) & = & -\sum_{\zeta=\e\,\eta} \li_3(\zeta, \eta)\\
\label{rel:23}
\sum_{\zeta\,\eta=\e} \li_2(\zeta, \eta) & = & - \sum_{\zeta\,\eta=\e} \li_3(\zeta, \eta)
\end{eqnarray}

We are now ready to start the proof of  Proposition~\ref{Prop:basis}.
It states that the lines $\li_j(\cdot,0), \li_j(\cdot,m-1)$ are superfluous in the sense that the remaining lines already generate the span of {\em all} lines, including these superfluous ones. 
In other words, Proposition \ref{Prop:basis} claims that these superfluous lines 
can be expressed as linear combinations of the remaining lines in $\NS(S)\otimes \Q$. 
To prove this, we work with the $6m\times 6m$-matrix $M$ whose entries are the coefficients of the superfluous lines in the relations (\ref{rel:H'}) and (\ref{rel:12})--(\ref{rel:23}).

The entries of the matrix $M$ are ordered as follows:

\begin{tabular}{lll}
columns & lines & $\li_1(0,0),\hdots,\li_1(m-1,0), \li_1(0,m-1),\hdots, \li_1(m-1,m-1)$,\\ 
&& $\li_2(0,0),\hdots,\li_3(m-1,m-1)$\\
rows & relations & (\ref{rel:H'}) for $\eta=\gamma^l, l=0,\hdots, m-1$ and $j=1,2,3$\\
&& (\ref{rel:12})--(\ref{rel:23}) for $\e=\gamma^i, i=0,\hdots,m-1$
\end{tabular}

That is to say, the matrix $M$ encodes the following system of relations on $\NS(S)$
\begin{eqnarray}\label{eq:M}
M\cdot \boldsymbol\li = \boldsymbol r
\end{eqnarray}
where the vector $\boldsymbol\li$ has entries the superfluous lines (ordered as above) and 
the right hand side vector $\boldsymbol r$ comprises the remaining terms of the chosen relation with the appropriate signs.

By the relations, all entries of $M$ are either $0$ or $1$. 
It will be convenient to write $M$ as a block matrix whose entries are 36 matrices of type $m\times m$.
In fact, the blocks arising from relation (\ref{rel:H'}) are just the identity Matrix $I$.
For the other relations, we need two permutation matrices of order $m$ which are transposes of each other:
\[
D=\begin{pmatrix} 0 & 1 & 0 & \hdots & & 0\\
0 & 0 & 1 & 0 & \hdots & 0\\
&& \hdots &&&\\
&&& \hdots &&\\
0 && \hdots && 0 & 1\\
1 & 0 && \hdots && 0
\end{pmatrix}, \;\;\;\;
B = D^{\rm t} = D^{-1}
\]
Then $M$ is given as follows: 
\[
M=\begin{pmatrix}
I & I & 0 & 0 & 0 & 0\\
0 & 0 & I & I & 0 & 0\\
0 & 0 & 0 & 0 & I & I\\
I & B & I & B & 0 & 0\\
I & D & 0 & 0 & I & B\\
0 & 0 & I & D & I & D
\end{pmatrix}
\]
We claim that there is a solution to the system of relations (\ref{eq:M}) in $\NS(S)\otimes\Q$.
If the matrix $M$ were invertible, then this would follow immediately.
However, $M$ is not invertible, 
so we have to find a way to circumvent this problem.

Recall that we are looking for a solution in $\NS(S)\otimes\Q$. 
Hence we can still modify any relation in $\NS(S)$ by adding
multiples of the relations (\ref{rel:H}) for any index $i$ and $\eta =1$ or $\eta=\gamma^{m-1}=\gamma^{-1}$. 
On the system of relations (\ref{eq:M}), this has the effect of adding a constant row to any of the six blocks 
associated to the invariants $i$ and $\eta$ of the chosen relation (\ref{rel:H}).
We will refer to this as adding constant rows. 
Of course, this modification changes the vector $\boldsymbol r$ on the right-hand side of (\ref{eq:M}) by adding a multiple of $H$, but we will not need to consider this expression at all. 

We will achieve a proof of Proposition \ref{Prop:basis}  by making the matrix $M$ invertible by adding constant rows.
First we shall simplify the matrix.
Note that elementary operations of linear algebra, if performed blockwise, are compatible with the modifications by adding constant rows.
This simplifies the problem of invertibility greatly:
\[
M=
\begin{pmatrix}
I & I & 0 & 0 & 0 & 0\\
0 & 0 & I & I & 0 & 0\\
0 & 0 & 0 & 0 & I & I\\
I & B & I & B & 0 & 0\\
I & D & 0 & 0 & I & B\\
0 & 0 & I & D & I & D
\end{pmatrix}
\to
\begin{pmatrix}
I & I & 0 & 0 & 0 & 0\\
0 & 0 & I & I & 0 & 0\\
0 & 0 & 0 & 0 & I & I\\
0 & B-I & 0 & B-I & 0 & 0\\
0 & D-I & 0 & 0 & 0 & B-I\\
0 & 0 & 0 & D-I & 0 & D-I
\end{pmatrix}
\]
\[
\to
\begin{pmatrix}
B-I & B-I & 0\\
D-I & 0 & B-I\\
0 & D-I & D-I
\end{pmatrix}
\to
\begin{pmatrix}
B-I & 0 & 0\\
D-I & 2I-B-D & B-I\\
0 & 0 & D-I
\end{pmatrix}
\]
To show that each of the superfluous lines can be expressed in terms of the other lines in $\NS(S)\otimes\Q$, it thus suffices to modify the following block matrices by adding constant rows such that they become invertible:
\[
B-I,\;\;\; D-I,\;\;\; 2I-B-D.
\]

\begin{Lemma}\label{Lem:ranknew}
\label{Lem:rank}
Let $U(r)$ denote the $m\times m$ matrix with entries $1$ in the $r$-th row and $0$ elsewhere. 
\begin{enumerate}[(i)]
\item
The determinants of $B-I+U(r)$ and $D-I+U(r)$ equal $(-1)^{m-1}m$ for any $r=1,\hdots,m$. 
\item
The determinant of $2I-B-D+U(2)$ equals $m^2$.
\end{enumerate}
\end{Lemma}

\begin{proof}
$(i)$
We calculate the determinants by computing all eigenvalues of the given matrices.
We claim that the eigenvalues are exactly
\begin{eqnarray}
\label{eq:ev}
\{\varepsilon-1; \varepsilon^m=1, \varepsilon\neq 1\}\cup \{1\}.
\end{eqnarray}
Then the determinant equals the product of the eigenvalues which can be written as
\[
\prod_{\varepsilon\neq 1} (\varepsilon-1)=\prod_{\varepsilon\neq 1} (\varepsilon-t)|_{t=1} = (-1)^{m-1} \left[\dfrac{t^m-1}{t-1}\right]_{t=1}=(-1)^{m-1} m.
\] 

To prove the claim about the eigenvalues, we exhibit simultaneous eigenvectors for all matrices $D, B, I, U(r)$.
This is easily accomplished by working with both multiplication from left and right.

For multiplication from the left, we have the common eigenvectors 
\[
\mathbf{v}_\varepsilon = (\varepsilon^i)_{0<i\leq m} \;\;\; \forall \; \varepsilon\in\mu_m\setminus\{1\}.
\]
These eigenvectors have eigenvalues $\varepsilon, \varepsilon^{-1}, 1, 0$, respectively.
Hence we obtain all eigenvalues from (\ref{eq:ev}) except for $1$.
The remaining eigenvalue is easily computed for
multiplication from the right.
Here we have the eigenvector
\[
\mathbf{v}_{1} = (1,\hdots,1)
\]
with eigenvalue 1 for each matrix $D,B,I,U(r)$. 
Thus the given matrices have the eigenvalue $1$.
This completes the proof of $(i)$.

For $(ii)$, note that
\[
B\cdot U(1) = U(2) \;\; \text{ and } \;\; U(1)\cdot D = U(1) \cdot U(1) = U(1).
\]
Together with the equality $DB = I$, this implies that 
\[
(B-I+U(1))\cdot (D-I+U(1)) = 2I-B-D + U(2).
\]
By $(i)$, this matrix has determinant $m^2$.
\end{proof}

By Lemma \ref{Lem:rank}
the matrix $M$ can be modified by adding constant rows to its blocks in such a way that it becomes invertible over $\Q$. 
Thus we can express all superfluous lines rationally in terms of the lines in $\mathcal{B}$. Since lines generate $\NS(S)$ rationally by Theorem~\ref{Thm:AS} and $\#\mathcal{B}=\rho(S)$, this completes the proof of Proposition \ref{Prop:basis}. 
\qed

\begin{Remark}\label{almostintegral}
\begin{enumerate}[(i)]
\item
The result of Proposition \ref{Prop:basis} stays valid in positive characteristic if the Picard number does not increase upon reduction (for instance for characteristics $p\equiv 1\mod m$).
\item
The method of proof does not require that $(m,6)=1$, but only that $m$ is odd.
For arbitrary odd degree $m$, we deduce that the lines in $\mathcal{B}$ generate the span of all lines $L$ rationally.
\item
For even degrees $m$, the matrix $M$ takes a different shape, as we cannot choose $\omega=-1$.
Hence the relations for $\alpha\in\D_m^1\cap\D_m^3$ change to
\[
w_1(\alpha) = - \omega^{2a_0}\,w_3(\alpha).
\]
Summing up as for odd $m$, we obtain
\[
\sum_{\zeta\,\eta=\e} \li_1(\zeta, \eta) = -\sum_{\omega^2\,\zeta=\e\,\eta} \li_3(\zeta, \eta)
\]
yielding a different relation matrix.
\end{enumerate}
\end{Remark}

\begin{Corollary}
\label{Cor:disc-B}
Let $m$ be any odd integer. Let $\Lambda \subset \NS(S)$ be the lattice generated by all lines and $\Lambda'$ the sublattice generated by those in $\mathcal{B}$  of Proposition \ref{Prop:basis}. 
Then the index $[\Lambda: \Lambda']$ is only divisible by primes dividing $m$.
In particular, $\Lambda'$ has discriminant dividing some power of $m$.
\end{Corollary}

\begin{proof}
The second claim follows from the first in conjunction with Corollary \ref{Cor:discr-m}.
For the first claim, it suffices to deduce from Lemma \ref{Lem:rank} that the matrix $M$ can be modified in such a way that it becomes invertible over $\Z[\frac 1m]$.
\end{proof}

\begin{Remark}\label{BtoLines} 
The modified matrices in Lemma \ref{Lem:rank} have determinant of absolute value $m$ or $m^2$.
There is no obvious way to make the matrix $M$ invertible over $\Z$.
Note, however, that we may still have $\Lambda' = \Lambda$ and even $\Lambda' = \NS(S)$, since the expression on the right-hand side of (\ref{eq:M}) might be divisible in $\NS(S)$.
In the cases of this paper with $(m,6)=1$, these equalities do
indeed hold. This  will be checked as part of the proof of Theorem~\ref{thm}.
\end{Remark}

For all odd degrees $m\leq 81$, we calculated the determinant of the intersection form of the lines in $\mathcal{B}$. In each case, the determinant turned out to be a perfect power of $m$, with exponent as conjectured in \cite{Sh-Jacobi}:
\begin{equation}\label{discB}
\det(\li.\li')_{\li,\li'\in\mathcal{B}} = m^{3(m-3)^2}.
\end{equation}

\section{Supersingular reduction technique}
\label{s:tech}

Consider the reduction of the complex Fermat surface $S$ mod $p$. Denote the resulting surface by $S_p$. Then $S_p$ is smooth for any $p\nmid m$. For any such $p$, reduction induces a specialisation embedding (see \cite[Proposition 3.6]{maulik}, and note that $\NS(S)$ and $\NS(S_p)$ are torsion-free)
\begin{eqnarray}\label{eq:NS-mod}
\NS(S) \hookrightarrow \NS(S_p).
\end{eqnarray}
We call a surface $X$ supersingular if its Picard number is maximal: $\rho(X)=b_2(X)$. For Fermat surfaces, we have the following result of Katsura and Shioda:

\begin{Theorem}[Katsura-Shioda {\cite{KS}}]\label{Thm:KS}
The reduction $S_p$ is supersingular if and only if there is some $r\in\N$ such that 
\[
p^r \equiv -1\mod m.
\]
\end{Theorem}

One advantage of working with supersingular surfaces is that we have good knowledge about the discriminant of their N\'eron-Severi groups. The following result is a generalisation of Artin's classification of supersingular K3 surfaces \cite{Artin}.

\begin{Theorem}[Ekedahl {\cite{teke}}, Sch\"utt--Schweizer {\cite{SS}}]
\label{Thm:NS}
Let $X$ be a smooth projective surface over a finite field $k$ of
characteristic $p$. Assume that $X$ is supersingular. 
Then
\[
|\disc(\Num(X))| =
p^{2\sigma}\;\;\;\;(\sigma\in{\N}_0).
\]
\end{Theorem}

The proof in \cite{SS} uses exactly the same techniques as Artin's original paper, mainly the Artin-Tate conjecture. The proof in \cite{teke} is based on cohomological results by Illusie and even allows to compute the (Artin) invariant $\sigma$.

We now explain the method by which we will prove Theorem~\ref{thm}. For this we recall the second betti number of $S$:
\[
b=b_2(S)=m^3-4m^2+6m-2.
\]
We shall also use the Lefschetz number $\lambda(S)=b_2(S)-\rho(S)$.

\textbf{Supersingular reduction technique}\\
{\sl
Fix the degree $m$. Let $p$ be a prime of supersingular reduction for $S$.

\begin{enumerate}
\item Compute a basis of $\NS(S)\otimes{\Q}$ consisting of lines $\li_j$. 

\item Let $N=\langle \li_j; j=1,\hdots,\rho\rangle \subseteq \NS(S)$. Compute disc$(N)$ in terms of the Gram matrix of the intersection numbers of the lines. Then disc$(N)=\nu^2\,\disc(\NS(S))$ where $\nu$ denotes the index of $N$ in $\NS(S)$.

\item Complement the reductions of the lines $\li_j (j=1,\hdots,\rho)$ by $\lambda(S)$ divisor classes $d_k$ on the supersingular reduction $S_p$ for a basis of $\NS(S_p)\otimes{\Q}$.

\item Let $N_p=\langle \li_j, d_k; \; j=1,\hdots,\rho; \; k=1,\hdots,b-\rho\rangle \subseteq \NS(S_p)$. Compute disc$(N_p)$.   
\end{enumerate}}

If $(m,6)=1$, then we will work with the rational basis $\mathcal{B}$ from Proposition~\ref{Prop:basis} in step 1. At the end of the previous section, we computed the discriminants of the lattice $N$ generated by these lines for several $m$. Recall that this discriminant was always a power of $m$ (and in general it is a divisor of some power of $m$ by Corollary \ref{Cor:disc-B}).

\begin{Criterion}\label{Crit}
Assume that the discriminants $N$ and $N_p$ have squarefree greatest common divisor.
Then $N=\NS(S)$ (i.e.~$\nu=1$).
\end{Criterion}

\begin{proof}
Let $D\in \NS(S)$. Consider the lattices
\[
N' = \langle  N, D\rangle ,\;\;\; N_p' = \langle  N_p, D\rangle.
\]
Let $r=[N':N]$, i.e.~$r$ is the minimal positive integer such that $rD\in N$, 
and we can write in $N$
\begin{eqnarray}\label{eq:rD}
rD = \sum a_i \li_i \;\;\; (a_i\in\Z).
\end{eqnarray}
We claim that this implies $r=[N_p':N_p]$. 
Assume on the contrary that there is a positive integer $s<r$ with $sD\in N_p$.
By assumption, we can write in $N_p$
\begin{eqnarray}\label{eq:sD}
sD = \sum b_i \li_i + \sum c_k d_k \;\;\; (b_i, c_k \in\Z).
\end{eqnarray} 
Necessarily there is some index $k$ with $d_k\neq 0$, 
for otherwise (\ref{eq:sD}) would be a relation in $N$, thus contradicting the minimality of $r$.
As not all $d_k$ are zero, the equations (\ref{eq:rD}) and (\ref{eq:sD}) combine to a non-trivial relation between the basis elements $\li_i, d_k$ of $N_p$.
This is impossible, hence the index of $N_p$ in $N_p'$ is $r$ as claimed.

We conclude that the lattices $N', N_p'$ have discriminants
\[
\disc(N') = \disc(N)/r^2,\;\;\; \disc(N_p') = \disc(N_p)/r^2.
\]
As the discriminants are integers, $r^2$ divides the greatest common divisor of the discriminants of $N$ and $N_p$.
By assumption, $r=1$ and hence $D\in N$.
\end{proof}

In sections \ref{s:4}--\ref{s:13}, we will apply the supersingular reduction technique to the Fermat surfaces of degree $4, 5, 7, 11$ and $13$. 
For a generalisation of Criterion~\ref{Crit}, 
one should note that the above proof does not actually require that $N_p$ has finite index in $\NS(S_p)$.
Hence we can also apply the same technique to sublattices of positive corank in $\NS(S_p)$ 
(which is computationally preferable as we can work with lattices of substantially smaller rank).
This approach will be extended in section~\ref{s:out} before we apply it to the degrees $m\geq 17$ in order to complete the proof of Theorem \ref{thm}.

\subsection{Additional lines mod $p$}
\label{s:add}

The supersingular reduction technique requires to complement the lines from characteristic zero by divisors which only appear after reduction modulo a supersingular prime $p$.
In this section, we will show how one can exhibit such divisors.
We concentrate on the case where the degree equals $q+1$ for some prime power $q=p^r$.
In general, this situation can be achieved by replacing the degree $m$ by a suitable multiple $mk$.
Then one can map down the divisors on the Fermat surface $\hat S_p$ of degree $mk$ to $S_p$ by the $k$-th power map
\begin{eqnarray*}
\hat S_p & \to & S_p\\
x_i & \mapsto & x_i^k.
\end{eqnarray*}

Throughout this section, we let $p$ be a prime, $r\in\N$ and $q=p^r$.
We fix the degree $m=q+1$ of the Fermat surface $S_p$ and perform our calculations over $\F_q$. 
In this situation, Tate and Thompson realised that the unitary group over $\F_{q^2}$ acts irreducibly on the primitive part of $H^2(S_p)$ (cf.~\cite{TT}). This provided the first proof for the if-part in Theorem \ref{Thm:KS}. 
In consequence, the images of any line on $S_p$ under the action of the unitary group generate $\NS(S_p)$ rationally together with the hyperplane section. 

In the sequel, we shall exhibit very specific lines for different choices of $m>3$.
In each case, we shall only give one line. 
Many further lines are obtained by applying the automorphisms of the surface to this line.
For our purposes, it will suffice to consider the images under the abelian group $\mu_m^4/\mu_m$ studied before.

\subsection{General $m$}
\label{ss:m}

Let $\alpha\in\F_q^*$ with $\alpha^2\neq -1$. Then consider the solutions $\beta\in\F_{q^2}$ of 
\begin{eqnarray}\label{eq:beta}
\beta^2=1+\alpha^2.
\end{eqnarray}
Since $m-2=q-1$, we have $\alpha^{m-2}=1$. 
As $\beta^2\in\F_q^*$, we also have
\[
\beta^{2\,(m-2)}=1.
\]
There are at least two $\alpha\in\F_q^*$ such that each solution $\beta$ of (\ref{eq:beta}) satisfies
\[
\beta^{m-2}=-1.
\]
For each such pair $(\alpha,\beta)$, we obtain the following line on $S_p$:
\[
\li_p=\{[\lambda, \alpha\lambda+\beta\mu, \beta\lambda+\alpha\mu,\mu];\;\; [\lambda,\mu]\in\PP^1\}.
\]

For many $m=q+1$, we can find simpler lines on $S_p$. We consider two cases:

\subsection{${m\equiv 2\mod 3}$}
\label{ss:2mod3}

If $m\equiv 2\mod 3$, i.e.~$q\equiv 1\mod 3$, then let $\alpha\in\F_q$ be a primitive third root of unity: $\alpha^2+\alpha+1=0$. Then $S_p$ contains the following line:
\[
\li_p = \{[\lambda, \alpha(\lambda+\alpha\mu), \alpha(\alpha\lambda-\mu),\mu];\;\; [\lambda,\mu]\in\PP^1\}.
\]

\subsection{${p=3}$}
\label{ss:p=3}

Let $p=3$. For any $q=p^r$ and $m=q+1$, $S_p$ contains the following line:
\[
\li_p = \{[\lambda, (\lambda+\mu), (\lambda-\mu),\mu];\;\; [\lambda,\mu]\in\PP^1\}.
\]

\subsection{Notation}
\label{s:act}

In the sequel, we shall always fix one line $\li_p$ as above. Then we let the subgroup $\mu_m^4/\mu_m$ of $\mbox{Aut}(S)$ act on $\li_p$. For convenience, we normalise the action of $\mu_m^4/\mu_m\cong \mu_m^3$ corresponding to the choice $\zeta_3=1$:
\[
g=(\zeta, \eta, \xi)\in\mu_m^3:\;\;\; [x_0,x_1,x_2,x_3] \mapsto [\zeta\,x_0,\eta\,x_1, \xi\,x_2,x_3].
\]
As before, we denote the resulting $m^3$ lines by
\[
\li_p(\zeta, \eta, \xi)=g(\li_p) \;\;\; \text{ or } \;\;\; \li_p(j,k,l) \;\; \text{ if } \zeta=\gamma^j, \eta=\gamma^k, \xi=\gamma^l.
\]
To identify the latter  lines, we shall always consider the reduction of the primitive root of unity $\gamma\in\mu_m$ that was used to enumerate the lines 
$\li_j(k,l)$ on $S$ in characteristic zero.

\begin{Remark}
\label{Rem:zero}
In the supersingular case, $V(\alpha)\subset H^2(S_p)$ is algebraic for any character $\alpha\in\A_m$. Given a line $\li_p$ as above, we can mimic the construction from section \ref{s:rat} to produce an eigendivisor with character $\alpha=(a_0,a_1,a_2,a_3)$:
\[
w_p(\alpha) = \sum_{\zeta,\eta,\xi} \zeta^{a_0}\,\eta^{a_1}\,\xi^{a_2}\,\li_p(\zeta,\eta,\xi).
\]
However, it is non-trivial to decide whether $w_p(\alpha)$ is non-zero in $\NS(S_p)$ (cf.~Remark~\ref{Rem:133}).
\end{Remark}

\section{Fermat surfaces of low degree}
\label{s:low}

In this section, we give a proof of Theorem \ref{thm} for degrees $m=4,5,7,11,13$ that is based on the supersingular reduction technique.
For $m=4$, this result has been known since the mid 70's.
We will review the historical development and give an alternative proof.
For $m>4$, the result is new.

\subsection{The Fermat quartic revisited}
\label{s:4}

In this section, we let $m=4$. Thus $S$ is a singular K3 surface (in the sense that $\rho(S)=20$, the maximum possible over $\C$). It was shown by Pjatecki\u\i -\v Sapiro and \v Safarevi\v c \cite{PSS} that $\NS(S)$ has discriminant $d=-16$ or $-64$. The latter is the case if the N\'eron-Severi group is generated by lines. 
Depending on a claim by Demjanenko, Pjatecki\u\i -\v Sapiro and \v Safarevi\v c    deduced  $d=-64$. 
However, Demjanenko's argument contained a mistake. A correction was given by Cassels in 1978 \cite{Cassels}.

In the meantime, Mizukami had investigated the following family of K3 surfaces:
\[
X_\lambda:\;\;\; \{x^4+y^4+z^4+w^4 = 2\,\lambda\,(x^2\,y^2+z^2\,w^2)\}\subset\PP^3.
\]
The following result was part of his Master's thesis in 1975 \cite{Mizukami}:

\begin{Proposition}[Mizukami]
Let $X_\lambda$ as above. Then $\rho(X_\lambda)\geq 19$, and 
\[
 \disc (\NS(X_\lambda)) = 
\begin{cases}
 -64, & \text{if } \lambda=0,\\
128, & \text{if } \rho(X_\lambda)=19.
\end{cases}
\]
\end{Proposition}

For the Fermat quartic, this result implied $d=-64$.
Thus it follows that lines generate $\NS(S)$ integrally (Proposition~\ref{Prop:4}).
An alternative proof can be based on another result about certain Kummer surfaces by Inose \cite{Inose}. 

Here we present an alternative argument using the supersingular reduction technique from section~\ref{s:tech} at the prime $p=3$. 
Note that by Theorem~\ref{Thm:KS} a prime $p$ is supersingular if and only if $p\equiv 3\mod 4$.
Since $m$ is even, the situation differs from the cases considered in section \ref{s:lines-basis}.
In particular, we cannot use $\omega=-1$;
instead we need $\omega$ with $\omega^4=-1$, so that we can use $\gamma = \omega^2$.

\begin{enumerate}
\item 
A rational basis $\mathcal{B}'$ of $\NS(S)$ can be expressed in terms of $\mathcal{B}$ as in Proposition \ref{Prop:basis} by switching $l\mapsto l-1$ and adding $\li_2(0,m-2)$: 
\[
\mathcal{B}' = \{\li_j(k,l); \li_j(k,l+1)\in\mathcal{B}\}\cup \{\li_2(0,m-2)\}.
\]

\item Let $N=\langle \li; \li\in \mathcal{B}'\rangle $. Then discr$(N)=-64$. 

\item On the supersingular reduction $S_3$, we have the additional line
\[
\li_3 = \{[\lambda, (\lambda+\mu), (\lambda-\mu),\mu];\;\; [\lambda,\mu]\in\PP^1\}
\]
from section \ref{ss:p=3}. Recall $\gamma$, the fixed square root of $-1$. Let
\[
\li_3'=  \{[\lambda, \gamma\,(\lambda+\mu), (\lambda-\mu),\mu];\;\; [\lambda,\mu]\in\PP^1\}.
\]       
Then we compute that the lines $\li\in \mathcal{B}'$ together with $\li_3, \li_3'$ constitute a rational basis $\mathcal{B}_3$ of $\NS(S_3)$:

\item Let $N_3=\langle \li; \; \li\in\mathcal{B}_3\rangle $. Then discr$(N_3)=-9$.  

\end{enumerate}
By Criterion~\ref{Crit}, we deduce that $N=\NS(S)$. In other words we have reproven the following result:

\begin{Proposition}[Mizukami, Inose]\label{Prop:4}
The complex Fermat quartic surface has N\'eron-Severi group generated by lines. Its discriminant is $-64$.
\end{Proposition}

%

The next result was first pointed out to the second author by Mizukami in the 1970's (unpublished report). Mizukami's proof was based on the computation of the intersection matrix for a suitable collection of lines on $S_3$.

\begin{Lemma}[Mizukami]
The reduction $S_3$ of the Fermat quartic mod $3$ has N\'eron-Severi group generated by lines over $\F_9$.
\end{Lemma}

\emph{Proof:} Since $S_3$ is a supersingular K3 surface, the exponent $\sigma$ from Theorem~\ref{Thm:NS} is the Artin invariant of $S_3$. By Artin's stratification \cite{Artin}, $\sigma\in\{1,\hdots,10\}$. Since the sublattice $N_3$ of $\NS(S_3)$ has discriminant $-9$, we deduce $N_3=\NS(S_3)$. \qed

\subsection{Fermat quintic}
\label{s:5}

In this section we shall prove Theorem~\ref{thm} for the complex Fermat quintic surface $S$. Note that $\rho(S)=37, b_2(S)=53$. 
It follows from Theorem~\ref{Thm:KS} that $p=2$ is a supersingular prime. We now apply the supersingular reduction technique from section~\ref{s:tech}.

\begin{enumerate}
\item
Take the rational basis $\mathcal{B}$ of $\NS(S)$ from Proposition~\ref{Prop:basis}.

\item 
Then $N=\langle \li; \li\in\mathcal{B}\rangle $ has discriminant $5^{12}$.
\end{enumerate}

On the supersingular reduction $S_2$ mod $2$, section \ref{ss:2mod3} gives 125 additional lines $\li_2(j,k,l)$ (plus their conjugates with respect to $\alpha\mapsto \alpha^2$). 
Here we write the third root of unity $\alpha$ in terms of a primitive fifth root of unity $\gamma$ as $\alpha=\gamma^3+\gamma^2+1$.

We express the 125 lines relative to $\gamma$ and $\alpha$ through one parameter $\nu=1,\hdots,125$ as $\li_p(j,k,l)=\li_p(\nu)$ where
\[
\nu=\nu(j,k,l)=25j+5k+l+1.
\]

\begin{enumerate}
\item[(3)]
Let $\mathcal{N}=\{32, 33, 34, 35, 36, 37, 38, 39, 44, 80, 81, 82, 83, 84, 93, 95\}$ and $\mathcal{B}_2=\{\li_p(\nu); \nu\in \mathcal{N}\}$.
Then $\mathcal{B}\cup\mathcal{B}_2$ constitutes a rational basis of $\NS(S_2)$.

\item[(4)] 
Let $N_2=\langle \li; \li\in\mathcal{B}\cup\mathcal{B}_2\rangle $. Then discr$(N_2)=2^{16}$.  
\end{enumerate}

By Criterion~\ref{Crit}, we deduce that $N=\NS(S)$ with discriminant $5^{12}$. In other words we have proven Theorem~\ref{thm} for the Fermat quintic surface.

By \cite[p.~12]{teke}, $\NS(S_p)$ has discriminant $p^{16}$ for all primes $p\equiv 2,3\mod 5$.
Hence we deduce

\begin{Lemma}
The N\'eron-Severi group  of the reduction of the Fermat quintic modulo $2$ is generated by lines over $\F_{16}$.
\end{Lemma}

\subsection{Fermat septic}
\label{s:7}

The Fermat septic surface $S$ has $\rho(S)=91, b_2(S)=187$. In characteristic zero, we have
\begin{enumerate}
\item
rational basis $\mathcal{B}$ of $\NS(S)$ from Proposition~\ref{Prop:basis},

\item 
lattice $N=\langle \li; \li\in\mathcal{B}\rangle $ of discriminant $7^{48}$.
\end{enumerate}

Since section \ref{s:add} only applies to $m=q+1$ for some prime power $q$, 
the Fermat septic $S$ does not admit any supersingular reduction with apparent additional lines. Instead we consider a suitable covering Fermat surface and push down the additional lines on a supersingular reduction.

Here we can work with the Fermat surface $\hat S$ of degree 14 and consider the reduction $\hat S_p$ mod $p=13$. In order to define a line mod $p$, we fix a primitive root $\gamma\in\mu_{7}$ as a zero of $x^2+5x+1$.
Let $\li_p$ denote the line from \ref{ss:m} for $\alpha=2, \beta=3\gamma+1$. 
Denote the push-down to $S$ by $D_p$. 
Then $D_p^2=-8$ by the adjunction formula.
The action of $\mu_7^4/\mu_7$ as in section \ref{s:act} gives divisors $D_p(j,k,l)$. We compute the following rational basis of $\NS(S_p)$:

\[
\mathcal{B}_p = \left\{ D_p(j,k,l);\;\; (j,k,l)\in I\right\}
\]
where 
\begin{eqnarray*}
I & = & I_1 \cup I_2\\
I_1 & = & \{(j,k,l);\;\; 0\leq j,k<m-1, 0<l<m-1\}\\
I_2 & = &  \{(j,0,0);\;\;0\leq j<m-1\} \cup \{(m-1,m-2,m-2)\}.
\end{eqnarray*}
The discriminant of the intersection form of the divisors in $\mathcal{B}_p$ is $2^{38}\,7^2\,13^{48}$. 

In order to combine the above divisors with the original lines from characteristic zero, we number them as follows:
\begin{eqnarray*}
I_1 \ni (j,k,l) & \mapsto & \nu(j,k,l) = 1+j+(m-1)\,k+(m-1)^2\,(l-1),\\
I_2 \ni (j,k,l) & \mapsto & \nu(j,k,l) = b_2(S) - (m-1) + j.
\end{eqnarray*}
With this notation, we can refer to $D_p(\nu)$ for $1\leq\nu\leq b_2(S)$. We then find a mixed basis using certain multiples of all $\nu$ in the range $1,\hdots,\lambda(S)$ modulo $b_2(S)$:

\begin{enumerate}
\item[(3)]
Let $\mathcal{N}=\{[31\,\nu \mod b_2(S)]; \;\; 1\leq\nu\leq \lambda(S)\}$ and $\mathcal{B}_p'=\{D_p(\nu); \nu\in\mathcal{N}\}$.
Then $\mathcal{B}\cup\mathcal{B}_p'$ constitutes a rational basis of $\NS(S_p)$.

\item[(4)] Let $N_p=\langle C; C\in\mathcal{B}\cup\mathcal{B}_p'\rangle $. Then discr$(N_p)=13^{40}$.  
\end{enumerate}

By Criterion~\ref{Crit}, we deduce that $N=\NS(S)$ with discriminant $7^{48}$. Thus we have proven Theorem~\ref{thm} for the Fermat septic surface. 

By \cite[p.~12]{teke}, the geometric genus $p_g(S)$ equals  the Artin invariant $\sigma$ of $S_p$ for all $p\equiv -1\mod m$ ($m$ being the degree of the Fermat surface $S$). 
For $m=7$ and $p=13$, the latter condition is fulfilled, and $p_g(S)=20$.
Hence we deduce $N_p=\NS(S_p)$.
In particular, it follows that $\NS(S_p)$ can be generated by divisors defined over $\F_{p^2}$.

\begin{Remark}
\label{Rem:133}
The choice $\alpha=1$ and $\beta=\sqrt{2}$ would yield another set of $m^3$ divisors on $S$. It is easily verified that the divisors from $\mathcal{B}_p$, even combined with the original lines from $\mathcal{B}$, only generate a sublattice of rank $133$ inside $\NS(S_p)$. 
This indicates that non-trivial linear combinations as in Remark~\ref{Rem:zero} might return zero for particular choices of $\alpha, \beta$.
\end{Remark}

\subsection{Fermat surface of degree $11$}
\label{s:11}

The Fermat surface $S$ of degree $m=11$ has $\rho(S)=271, b_2(S)=911$. In characteristic zero, we have
\begin{enumerate}
\item
rational basis $\mathcal{B}$ of $\NS(S)$ from Proposition~\ref{Prop:basis},

\item 
lattice $N=\langle \li; \li\in\mathcal{B}\rangle $ of discriminant $11^{192}$.
\end{enumerate}

Consider the supersingular reduction $S_p$ mod $p=2$.
In order to exhibit additional divisors on $S_p$, we consider the Fermat surface $\hat S$ of degree 33.
The covering map $\hat S\to S$ has degree $27$. By section \ref{s:add}, the reduction $\hat S_p$ admits many additional lines. These will be pushed down to $S_p$.

The primitive roots $\gamma\in\mu_{m}$ are given as zeroes of the irreducible polynomial $(x^m-1)/(x-1)$. Fix such a~$\gamma\in\F_{p^{10}}$. 
Let $\li_p$ denote the line from \ref{ss:m} for 
\[
\alpha=\gamma^8+\gamma^7+\gamma^6+\gamma^5+\gamma^4+\gamma^3,\;\;\; \beta=\alpha+1.
\]
Denote the push-down to $S$ by $D_p$. 
By the adjunction formula, as mentioned in section \ref{s:prelim}, we have $D_p^2=-23$. 
The action of $\mu_m^4/\mu_m$ as in section \ref{s:act} gives divisors $D_p(j,k,l)$. We compute the same rational basis $\mathcal{B}_p=\mathcal{B}_p(m)$ of $\NS(S_p)$ as in section \ref{s:7}. 
The lattice generated by the divisors in $\mathcal{B}_p$ has discriminant
\[
2^{1200}\,3^2\,11^2\,23^{64}\,43^{24}\,67^8\,131^{16}\,197^4\,307^8\,331^8\,463^{12}\,593^8\,3541^8.
\] 
With $m$ and $p$ replaced, we employ the same numbering of $D_p(\nu)$ for $1\leq\nu\leq b_2(S)$ as in the previous section. 
As before we determine a mixed basis by using appropriate multiples of all $\nu$ in the range $1,\hdots,\lambda(S)$ modulo $b_2(S)$:

\begin{enumerate}
\item[(3)]
Let $\mathcal{N}=\{[
253\,\nu \mod b_2(S)]; \;\; 1\leq\nu\leq \lambda(S)\}$ and $\mathcal{B}_p'=\{D_p(\nu); \nu\in\mathcal{N}\}$.
Then $\mathcal{B}\cup\mathcal{B}_p'$ constitutes a rational basis of $\NS(S_p)$.

\item[(4)] Let $N_p=\langle C; C\in\mathcal{B}\cup\mathcal{B}_p'\rangle $. Then $N_p$ has discriminant
\[
2^{1202}\,5^4\, 7^4\,23^{48}\,43^{16}\,131^{16}\,439^2.
\]  
\end{enumerate}
%
%
By Criterion~\ref{Crit}, we deduce that $N=\NS(S)$ with discriminant $11^{192}$. 
This completes the proof of Theorem~\ref{thm} for the Fermat surface of degree $11$.

\subsection{Fermat surface of degree $13$}
\label{s:13}

The Fermat surface $S$ of degree $m=13$ has $\rho(S)=397, b_2(S)=1597$. In characteristic zero, we have
\begin{enumerate}
\item
rational basis $\mathcal{B}$ of $\NS(S)$ from Proposition~\ref{Prop:basis},

\item 
lattice $N=\langle \li; \li\in\mathcal{B}\rangle $ of discriminant $13^{300}$.
\end{enumerate}

Consider the supersingular reduction $S_p$ mod $p=5$.
In order to derive additional divisors on $S_p$, we consider  the Fermat surface $\hat S$ of degree 26
which is a degree 8-covering of $S$. The reduction $\hat S_p$ admits many additional lines by section \ref{s:add}. 

Here, we fix a primitive root $\gamma\in\mu_{m}$ as a zero of $x^4+2\,x^3+x^2+2\,x+1$.
Let $\li_p$ denote the line from \ref{ss:m} for 
\[
\alpha = 2\gamma^3 + 2\gamma^2 + \gamma, \qquad
\beta=-\gamma^2  -\gamma + 3.
\]
Denote the 
push-down to $S$ by $D_p$. The action of $\mu_m^4/\mu_m$ as in section \ref{s:act} gives divisors $D_p(j,k,l)$. We compute the same rational basis $\mathcal{B}_p=\mathcal{B}_p(m)$ of $\NS(S_p)$ as in section \ref{s:7} and \ref{s:11}. The determinant of the intersection form of the divisors in $\mathcal{B}_p$ is 
\[
2^{26}\,3^{192}\,5^{912}\,13^2\,53^{24}\,79^{24}\,103^{32}\,181^8\,233^8\,313^8\,677^{16}\,883^4\,2003^8\,2729^8\,3847^8.
\] 
Employ the same numbering of $D_p(\nu)$ for $1\leq\nu\leq b_2(S)$. Again we find a mixed basis using appropriate multiples of all $\nu$ in the range $1,\hdots,\lambda(S)$ modulo $b_2(S)$:

\begin{enumerate}
\item[(3)]
Let $\mathcal{N}=\{[5\,\nu \mod b_2(S)]; \;\; 1\leq\nu\leq \lambda(S)\}$ and $\mathcal{B}_p'=\{D_p(\nu); \nu\in\mathcal{N}\}$.
Then $\mathcal{B}\cup\mathcal{B}_p'$ constitutes a rational basis of $\NS(S_p)$.

\item[(4)] Let $N_p=\langle C; C\in\mathcal{B}\cup\mathcal{B}_p'\rangle $. Then $N_p$ has discriminant
\[
2^4\, 3^{144}\, 5^{912}\, 53^{16}\, 103^{32}\, 677^{16}\, 1151^2\, 40627^2\, 42702482453593^2\, 247634616308749^2.
\]  
\end{enumerate}

By Criterion~\ref{Crit}, we deduce that $N=\NS(S)$ with discriminant $13^{300}$. This completes the proof of Theorem~\ref{thm} for the Fermat surface of degree $13$.

\section{Generalisations and extensions}
\label{s:out}

For Fermat surfaces of degrees up to $m=13$, we exhibited an explicit rational basis of $\NS(S_p)$ for some supersingular prime $p$, thus enabling us to apply the supersingular reduction technique.
This approach has two advantages: 
first we can double-check the compatibility with the discriminant of $\NS(S_p)$ from Theorem \ref{Thm:NS}; 
secondly we obtained additional information on generators of $\NS(S_p)$ in some cases.

For higher degrees, however, the matrices get too large for an explicit computation of the determinant.
In this section we develop an extension of Criterion~\ref{Crit}. 
This will allow us to treat much higher degrees and eventually give a full proof of Theorem \ref{thm}.
First we rephrase the old criterion in a more general setting. 

\begin{Lemma}\label{abeliangroups}
Suppose
$$
\xymatrix{
M \ar[r]^\varphi \ar[d]_\psi &L \ar[d] ^\chi\\
M' \ar[r]_{\varrho} & L'
}
$$
is a commutative diagram of homomorphisms of abelian groups with $\chi$ and $\varrho$ injective. 
Suppose that $L/\varphi(M)$ is torsion-free and that $M'/\psi(M)$ is torsion. Then $\varrho$ induces an injective homomorphism $M'/\psi(M) \to L'/\chi(L)$. If the group $L'/\chi(L)$ is finite, then the index $[M':\psi(M)]$ divides the index $[L':\chi(L)]$. 
\end{Lemma}
\begin{proof}
Set $\sigma = \varrho \circ \psi = \chi \circ \varphi$. As $\chi$ is injective, it induces an injection $\chi\colon L/\varphi(M) \to L'/\sigma(M)$. The quotient $(\chi(L) \cap \varrho(M')) / \sigma (M)$ is contained in $\chi(L/\varphi(M))$, which by injectivity of $\chi$ is torsion-free. The same quotient is also contained in $\varrho(M'/\psi(M))$, which is torsion. We conclude that the quotient is trivial, i.e., $\chi(L) \cap \varrho(M')  =  \sigma (M)$.
The kernel of the map $M' \to L'/\chi(L)$ induced by $\varrho$ is 
\[
\varrho^{-1}(\chi(L))= \varrho^{-1}(\chi(L) \cap \varrho(M')) = \varrho^{-1}(\sigma(M)) = \psi(M),
\] 
where the last equality follows from the injectivity of $\varrho$. The first statement of the lemma follows.
Assuming finiteness of $L'/\chi(L)$, the divisibility of indices follows immediately. 
\end{proof}

Recall that we only consider integral non-degenerate lattices.
The following proposition gives a method to show that a given lattice $M$ equals an a priori unknown superlattice $M'$ that contains $M$ as a sublattice of finite index.

\begin{Proposition}\label{lattices}
Suppose $\varrho \colon M' \to L'$ is an injective homomorphism of lattices. Let $M$ be a finite-index sublattice of $M'$ and $L$ a sublattice of $L'$ that contains $\varrho(M)$ primitively.  If the greatest common divisor $(\disc (M),\disc (L))$ is squarefree, then $M$ equals $M'$. 
\end{Proposition}
\begin{proof}
Let $L''$ be the saturation of $L$ in $L'$, i.e., $L'' = L_{{\Q}} \cap L'$, where the intersection is taken inside $L'_{{\Q}}$. From $M_{\Q}= M'_{\Q}$ we find 
\[
\varrho(M') \subset \varrho(M'_{\Q}) = \varrho(M_{\Q}) \subset L_{\Q}
\]
and conclude $\varrho(M') \subset L''$. After replacing $L'$ by $L''$, we may assume that $L$ has finite index in $L'$. By Lemma \ref{abeliangroups}, with $\varphi = \varrho$ and $\psi$ and $\chi$ being inclusions, we find that $[M':M]$ divides $[L':L]$.
From $\disc (M) = [M': M]^2 \disc (M')$ we conclude that $[M':M]^2$ divides $\disc (M)$ and similarly $[L':L]^2$ divides $\disc (L)$. Therefore $[M':M]^2$ divides $(\disc (M), \disc (L))$. If $(\disc (M), \disc (L))$ is squarefree, then it follows that $M$ equals $M'$.
\end{proof}

Criterion~\ref{Crit} is exactly Proposition \ref{lattices} applied to $M' = \NS(S)$ and $L' = \NS(S_p)$; 
the primitivity was ensured by complementing a basis of $M'_{\Q}$ to a basis of $L'_{\Q}$ (cf.~the proof of Criterion \ref{Crit}).
As mentioned at the end of section \ref{s:tech}, the sublattice in $L'$ does not need to have finite index in $L'$. 
In practice, Proposition \ref{lattices} will often be applied when we have $(\disc (M), \disc (L)) =1$. 
Suppose $\varrho \colon M' \to L'$ is an injective homomorphism of lattices whose discriminants we do not know. Assume we have a finite-index sublattice $M$ of $M'$ with discriminant $\Delta = \disc (M)$ that we do know and we wish to show that $M$ equals $M'$.
By Proposition \ref{lattices} it suffices to find a sublattice $L$ of $L'$ that contains $\varrho(M)$ primitively with $(\Delta, \disc (L))=1$ or more generally squarefree greatest common divisor.

\subsection{Alternative approach}

In the previous section, we suggested to use an intermediate lattice $\Lambda\subset L\subset\NS(S_p)$ for the supersingular reduction technique.
While this does decrease the size of the matrices considered, we still had to compute their determinants which may be infeasible.
Instead we shall pursue an alternative approach that decreases the size of the matrix drastically and has further computational advantages.
Before an abstract treatment of the method, we sketch the general idea for the Fermat surfaces.

Consider the Fermat surface $S$ of degree $m$ with $(m,6)=1$.
Let $\Lambda$ denote the sublattice of $\NS(S)$ generated by the lines in $\mathcal B$ as in Proposition \ref{Prop:basis}.
Suppose that $\Lambda\neq\NS(S)$, so there is a prime $\ell$ and a divisor $D_0\in\Lambda$ that is $\ell$-divisible in $\NS(S)$, but not in $\Lambda$.
Clearly this implies $\ell\mid(D_0.C)$ for any curve $C$ on $S$ -- and on $S_p$ for any prime $p$ of good reduction.

Now let $\mathcal{C}$ denote any finite subset of Div$(S)$ or Div$(S_p)$.
Then we build the matrix of intersection numbers
\[
Q = (D.C)_{D\in\mathcal B, C\in\mathcal C}.
\]
This matrix has integer entries, so we can also consider it over $\F_\ell$.

\textbf{Claim:}
The rank of $Q$ over $\F_\ell$ does not exceed $\#\mathcal B-1=\rho(S)-1$.

\begin{proof}
To see this, consider the map 
$$
\varphi \colon \Lambda_{\F_\ell} \to \Hom(\F_\ell^\mathcal{C} , \F_\ell)
$$
that sends $D \in \Lambda_{\F_\ell}$ to the map that sends $C \in \mathcal{C}$ to $(C\cdot D \mod \ell)$ (and is extended linearly to $\F_\ell^\mathcal{C}$). Then multiplication by $(Q \mod \ell)$ from the right describes the linear map $\varphi$ with respect to the basis $\mathcal{B}$ of $\Lambda_{\F_\ell}$ and the basis of $\Hom(\F_\ell^\mathcal{C} , \F_\ell)$ that is dual to $\mathcal{C}$. Since $D_0$ is not $\ell$-divisible
in $\Lambda$, its image in $\Lambda_{\F_\ell}$ is nontrivial. From $\varphi(D_0) = 0$ we conclude that $\varphi$ is not injective, so $Q$ does not have maximal rank over $\F_\ell$.

Alternatively, pick a basis containing the primitive closure $D'$ of $D$ in $\Lambda$.
Since $D'$ is still $\ell$-divisible in $\NS(S)$, all entries in the row of $Q$ corresponding to $D'$ are zero mod $\ell$.
Hence the rank of $Q$ over $\F_\ell$ cannot exceed $\#\mathcal B-1$.
\end{proof}

In order to show that $\Lambda=\NS(S)$, we find a suitable set $\mathcal C$ of divisors on $S$ or any good reduction $S_p$ such that the matrix $Q$ has maximal rank $\rho(S)$ over $\F_\ell$.
Since the index of $\Lambda$ in $\NS(S)$ divides an $m$-power by Corollary \ref{Cor:disc-B}, 
it suffices to carry out the above procedure for all prime divisors $\ell\mid m$.
This approach has several computational advantages:
\begin{enumerate}
\item
We can work with a relatively small matrix $Q$ of size $\rho(S)\times \#\mathcal C$.
\item
We can work with the matrix $Q$ mod $\ell$.
\item
The elements of $\mathcal C$ do not have to be independent in $\NS(S_p)$.
\item
We can add divisors to $\mathcal C$ successively until the kernel of multiplication by $Q$ (from the right) on $\F_\ell^{\rho(S)}$ is zero.
\end{enumerate}

We shall now give an abstract formulation of this approach.
In \ref{ss:Fermat}, we will apply the method to Fermat surfaces of degrees up to $m=97$ to complete the proof of Theorem \ref{thm}.

\subsection{Abstract formulation}

Suppose for this paragraph that the conditions of Proposition \ref{lattices} are met, 
so $L$ is a lattice containing $\varrho(M)$ primitively.
Let $\ell$ be a prime divisor of $\Delta$, the discriminant of $M$. 
The quotient $L/\varrho(M)$ is free, and it follows that the induced map $M_{\F_\ell} \to L_{\F_\ell}$ is injective. Since $\ell$ does not divide $\disc (L)$, the pairing $L_{\F_\ell} \times L_{\F_\ell} \to \F_\ell$ is nondegenerate in the sense that the induced map $L_{\F_\ell} \to L_{\F_\ell}^*$ is injective.  In particular, the restriction $M_{\F_\ell} \to L_{\F_\ell}^*$ is injective.  We will see that this is in fact sufficient to conclude $M=M'$. 

\begin{Proposition}\label{easy}
Let $\varrho \colon M \to L'$ be an injective homomorphism of lattices.
Suppose that for every prime $\ell$ dividing $\disc (M)$, there is a sublattice $L(\ell)$ of $L'$ containing $\varrho(M)$ such that the composition $M_{\F_\ell} \to L(\ell)_{\F_\ell}^*$ of the reduction $M_{\F_\ell} \to L(\ell)_{\F_\ell}$ of $\varrho$ with the map $L(\ell)_{\F_\ell} \to L(\ell)_{\F_\ell}^*$ induced by the pairing on $L(\ell)$, is injective. Then $\varrho(M)$ is primitively contained in $L'$. 
\end{Proposition}

\begin{proof}
Let $M'$ denote the saturation $\varrho(M)_{\Q} \cap L'$ of $\varrho(M)$ in $L'$, where the intersection is taken in $L'_{\Q}$.  Then the inclusion $M' \to L'$ induces an isomorphism 
$$
M'/\varrho(M) \to \big(L'/\varrho(M)\big)_\tors. 
$$
Let $\ell$ be a prime with $\ell\nmid \disc (M)$. From $[M' : \varrho(M)] \mid \disc (\varrho(M)) = \disc (M)$ we find 
$$
\ell \nmid [M' : \varrho(M)] = \#M'/\varrho(M) =\#\big(L'/\varrho(M)\big)_\tors,
$$
so the quotient $L'/\varrho(M)$ has no nontrivial $\ell$-torsion. 
Now let $\ell$ be a prime with $\ell\mid \disc (M)$ and consider the composition 
$$
M_{\F_\ell} \xrightarrow{\varrho_\ell} L(\ell)_{\F_\ell} \to L'_{\F_\ell} \to {L^{\prime *}_{\F_\ell}} \to L(\ell)_{\F_\ell}^*.
$$
Here $\varrho_\ell$ is the reduction of $\varrho$ mentioned in the proposition, the second map is the reduction of the inclusion $L(\ell) \subset L'$, the third is induced by the pairing on $L'$, and the last is the dual of the second. Then the composition of the last three maps is induced by the pairing on $L(\ell)$, so the full composition is injective by assumption.  This implies that the composition 
\[
\tau \colon M_{\F_\ell} \to L'_{\F_\ell}
\]
of the first two maps is injective. Suppose $y \in L'/\varrho(M)$ satisfies $\ell y=0$. Let $x \in L'$ be a lift of $y$, so that there is an $m \in M$ with $\varrho(m) = \ell x$. The reduction $\overline{m} \in M_{\F_\ell}$ satisfies $\tau(\overline{m}) = 0$, so by injectivity of $\tau$, we obtain $\bar m=0$, i.e.~there is an $m' \in M$ with $\ell m' = m$. Then we have 
\[
\ell\varrho(m') = \varrho(m) = \ell x,\;\; \text{ so }\; \ell(\varrho(m') -x) =0.
\]
As $L'$ is torsion-free, we conclude $\varrho(m') =x$ and thus $y=0$. We deduce that again $L'/\varrho(M)$ has no nontrivial $\ell$-torsion, and therefore that $L'/\varrho(M)$ is torsion-free, i.e., $\varrho(M)$ is contained primitively in $L'$.
\end{proof}

\begin{Corollary}\label{usefulcor}
Suppose $\varrho \colon M' \to L'$ is an injective homomorphism of lattices. 
Let $M$ be a finite-index sublattice of $M'$.
Suppose that for each prime $\ell$ dividing $\disc (M)$, there is a sublattice $L(\ell)$ of $L'$ containing $\varrho(M)$ such that the induced map $M_{\F_\ell} \to L(\ell)_{\F_\ell}^*$ is injective. 
Then $M$ equals $M'$.
\end{Corollary}
\begin{proof}
We have inclusions $\varrho(M) \subset \varrho(M') \subset L'$.  
By Proposition \ref{easy}, the lattice $\varrho(M)$ is primitively contained in $L'$, so also primitively in $\varrho(M')$. 
As $\varrho(M)$ has finite index in $\varrho(M')$, we find $\varrho(M) = \varrho(M')$
and thus $M=M'$ by injectivity of $\varrho$. 
\end{proof}

Corollary \ref{usefulcor}  is weaker than Proposition \ref{lattices} in the sense that it implicitly assumes that the map $M_{\F_\ell} \to L_{\F_\ell}^{\prime *}$ is injective.
For instance,
Corollary \ref{usefulcor} cannot be applied  in the case $M=M'=L' = \langle e_\ell\rangle$, where $\langle e_\ell\rangle$ denotes a one-dimensional lattice whose generator $e_\ell$ has norm $\ell$ for some prime number $\ell$; Proposition \ref{lattices} does apply, as $\disc (M) = \ell$ is squarefree. 

However, Corollary \ref{usefulcor} has several advantages over Proposition \ref{lattices}, especially computationally. First of all, we only need to know the pairing between elements in a basis $A$ for $M$ and those in a set $B$ of generators for $L=L(\ell)$, as opposed to the pairing among all elements of $B$, which saves a lot of work when the rank of $L$ is much larger than that of $M$. Furthermore, we do not need to compute the discriminant of the larger lattice $L$. This also means that we do not even need to find a basis among the elements of $B$. Also, all computations can be done over $\F_\ell$ instead of $\Z$, which for finding (large) ranks makes quite a difference.
Finally, Proposition \ref{easy} and Corollary \ref{usefulcor} can easily be modified in such a way that it is possible to work with different lattices $L'$ and embeddings $\varrho: M\to L'$ for each prime $\ell\mid\disc(M)$.
In the framework of the supersingular reduction technique, one could then take different supersingular primes of the Fermat surface $S$ for each prime divisor $\ell$ of the degree $m$.

\subsection{Application to surfaces}
\label{ss:app}

Now suppose $\mathcal X$ is a nice 
surface over ${\Z}[1/N]$ (so smooth, projective, and every geometric fiber is integral) for some integer $N$ and denote $X=\mathcal X\otimes \bar\Q$.
Let $p\nmid N$ be a prime, so that $p$ is a prime of good reduction of $\mathcal X$ and denote
$X_p = \mathcal X\otimes \bar\F_p$. 
Then there is an injective homomorphism 
$$
\NS(X)/(\ptors) \hookrightarrow \NS(X_p)/(\ptors)
$$
of lattices (see \cite[Proposition 3.6]{maulik}). 
We can therefore apply Proposition \ref{lattices} or Corollary \ref{usefulcor} with 
\[
M' = \NS(X)/\tors \cong \Num(X)\;\;\text{ and }\;\;L' = \NS(X_p)/\tors \cong \Num(X_p),
\]
while $M$ is a finite-index sublattice of $M'$.
 This means, that if a priori we do not yet know the lattice $\Num(X)$, but we do know its rank $\rho = \rk \Num(X) = \rk \NS(X)$ and a sublattice $M \subset \Num(X)$ of rank $\rho$, then this gives a method to prove that $\Num(X)$ equals $M$; it suffices to find a lattice $L$ as in Proposition \ref{lattices} (as we have done in the previous sections) or lattices $L(\ell)$ as in Corollary \ref{usefulcor}. Note that $L$ and $L(\ell)$ do not need to have the same rank as $L'$. If they do have the same rank, and thus finite index in $\Num(X_p)$, then this may also give extra information about $\Num(X_p)$, such as an upper bound for its discriminant.

\subsection{Fermat surfaces}
\label{ss:Fermat}

In order to complete the proof of Theorem \ref{thm}, we continue to consider the Fermat surfaces $\mathcal S=\mathcal S_m \subset \P^3$ 
over ${\Z}[\frac 1m]$ given by 
\[
x^m+y^m+z^m+w^m=0
\]
for any integer $m >4$ with $\gcd(m,6)=1$. 
As in \ref{ss:app}, we let $S=\mathcal S\otimes\bar\Q$ and $S_p=\mathcal S\otimes\bar\F_p$ for any prime $p\nmid m$.
Sometimes we will also indicate the degree $m$ as a subscript and write $S_m$ and $(S_m)_p$, but whenever the degree is clear, it will be omitted.
Then as before, $\NS(S)$ and $\NS(S_p)$ 
are torsion-free for any prime $p\nmid m$ (see section \ref{s:prelim}).

The following table contains for each $m$ with $4 < m < 100$ 
and $\gcd(m,6)=1$ an integer $r$ such that $q=rm-1$ is a prime power, namely $q = p^n$ with $p$ prime, a prime $\ell \mid m$, an irreducible polynomial of 
degree $2n$ over $\F_p$, and one or two pairs $(\alpha,\beta)\in \F_q \times \F_{q^2}$. A dash indicates the same as the entry above it.

\hspace{-2cm}
\begin{minipage}{6in}
$$
\begin{array}{|cccccc|}
\hline
m & r & p^n  &\ell & f & (\alpha,\beta)\\
\hline \hline
5   & 1 & 2^2  & 5 & (x^5-1)/(x-1) & (\gamma^3+\gamma^2+1,\alpha+1)  \\
\hline
7   & 2 & 13  & 7 & x^2+3x+1 & (11, 11\gamma+10)  \\
\hline
11 & 3 & 2^5  & 11 & (x^{11}-1)/(x-1) & (\gamma^9+\gamma^8+\gamma^3+\gamma^2+1, \alpha+1)  \\
\hline
13 & 2 & 5^2  & 13 & x^4+x^3-x^2+x+1 & (-\gamma^3-\gamma^2+2\gamma+1, 3\gamma^3+\gamma^2+3\gamma-1) \\
\hline
17 & 1 & 2^4  & 17 & x^8+x^5+x^4+x^3+1 & (\gamma^7+\gamma^5+\gamma^4+1,\alpha+1)  \\
\hline
19 & 2 & 37  & 19 & x^2+3x+1 &(13,5\gamma+26)  \\
\hline
23 & 6 & 137  & 23 & x^2+11x+1 & (67,91\gamma+21)  \\
\hline
25 & 2 & 7^2  & 5 &  x^4+2x^3+4x^2+2x+1 & (\gamma^3+2\gamma^2+3\gamma,5\gamma^3+\gamma^2-1) \\
\hline
29 & 6 & 173  & 29 & x^2+18x+1 & (137,127\gamma+105)  \\
\hline
31 & 2 & 61  & 31 & x^2+5x+1 & (-3,11\gamma-3)  \\
\hline
35 & 4 & 139  & 5 & x^2+4x+1 & (-15,86\gamma+33) \\
- & - &  - & 7 & -& - \\
\hline
37 & 2 & 73  & 37 & x^2+3x+1 & (31,5\gamma+44)  \\
\hline
41 & 2 & 3^4  & 41 & x^8+x^6+x^5-x^4+ \hfill& (\gamma^7+\gamma^6+2\gamma^4+\gamma^2+2, \gamma^7+2\gamma^6+2\gamma^3+\gamma+1) \\
      &    &            &      &  \hfill +x^3+x^2+1 & \\
\hline
43 & 3 & 2^7  & 43 &  x^{14}+x^{11}+x^{10}+x^9+x^8+ \hfill& (\gamma^{12}+\gamma^{11}+\gamma^9+\gamma^8+\gamma^6+\gamma^5+\hfill\\
      &     &     &      & +x^7+x^6+x^5+x^4+x^3+1 & \hfill +\gamma^4+\gamma^3+\gamma^2, \alpha+1) \\
\hline
47 & 6 & 281 & 47 & x^2+10x+1 & (-18,158\gamma+228)  \\
\hline
49 & 2 & 97  & 7 & x^2+3x+1 & (-6,7\gamma+59) \\
\hline
53 & 4 & 211  & 53 & x^2+4x+1 & (34,33\gamma+66)  \\
\hline
55 & 2 & 109  & 5 &  x^2+6x+1 & (72,12\gamma+36)\\
- & - & - & 11 & - & (53,18\gamma+54),(73,51\gamma+44) \\
\hline
59 & 6 & 353  & 59 &  x^2+3x+1 & (-28,236\gamma+1) \\
\hline
61 & 2 & 11^2  & 61 & x^4+x^3+3x^2+x+1& (\gamma^3+\gamma^2+2\gamma+8,-\gamma^3+6\gamma^2+3\gamma+4)\\
\hline
65 & 1& 2^6 &5 & x^{12}+x^8+x^7+x^6+ \hfill& (\gamma^{11}+\gamma^9+\gamma^7+\gamma^6+\gamma^3+\gamma^2,\alpha+1), \\
      &  & &    & \hfill+x^5+x^4+1 & (\gamma^9+\gamma^5+\gamma^4+\gamma^2+\gamma,\alpha+1) \\
      -  & - & -& 13 &  - & (\gamma^{10}+\gamma^9+\gamma^7+\gamma^6+\gamma^5+\gamma+1,\alpha+1)\\
\hline
67 & 6 & 401 & 67  &x^2+24x+1 & (222,229\gamma+342)\\
\hline
71 & 4 & 283 &  71 & x^2+4x+1 &(-39,160\gamma+37 ) \\
\hline
73 & 10 & 3^6 & 73 & x^{12}+x^{10}+x^7-x^6+ \hfill& (-\gamma^{11}+\gamma^{10}-\gamma^8-\gamma^5+\gamma^4+\gamma^3+\gamma+1,\hfill\\
      &           &   &     & \hfill +x^5+x^2+1 & \hfill \gamma^{11}+\gamma^9+\gamma^8+\gamma^7+\gamma^6+\gamma^5+\gamma^4-\gamma^3+\gamma^2+\gamma) \\
\hline
77 & 4 & 307  & 7& x^2+4x+1 &(29,136\gamma-35),(-73,61\gamma+122) \\
- & -& -& 11&-& (197,-51\gamma+205),(91,-10\gamma-20)  \\
\hline
79 & 2& 157& 79 & x^2+3x+1& (-5,127\gamma+112) \\
\hline
83 & 4& 331& 83 & x^2+4x+1& (163,19\gamma+38) \\
\hline
85 & 2 & 13^2& 5 &x^4+x^3+4x^2+x+1&(8\gamma^3 + 8\gamma^2 -2\gamma-1,7\gamma^3-\gamma^2+\gamma+4), \\
     &    &        &   &   &  (6\gamma^3+6\gamma^2+5\gamma,10\gamma^3-\gamma^2+1) \\
- & -& -& 17 & -& -\\
\hline
89 & 16 & 1423 & 89& x^2+14x+1 & (536,184\gamma+49) \\
\hline
91 & 2& 181 & 7 & x^2+5x+1&(145,139\gamma+76),(80,109\gamma+1) \\
 -& -&- &  13 &- & - \\
\hline
95 & 4 & 379 &  5 & x^2+59x+1& (35,243\gamma+157),(200,\gamma+219)\\
- & -&- &  19& -& (45,89\gamma+162),(26,8\gamma+236) \\
\hline
97 & 2 & 193  & 97 & x^2+3x+1 & (6,50\gamma+75)\\
\hline
\end{array}
$$
\end{minipage}

Before elaborating on how the contents of this table was computed, let us explain its meaning.
Suppose $m,r,p^n=q=rm-1$, $\ell,f$, and 
the $s$ pairs $(\alpha_1,\beta_1),\ldots,(\alpha_s,\beta_s)$ are the elements contained in one row of the table. 
\begin{itemize}
\item
Let $\gamma$ denote a root of $f$ in $\F_p[x]/(f) \cong \F_{q^2}$. 
By choice of $f$, $\gamma$ is a primitive $m$-th root of unity. 
\item
Let $M \subset \NS(S)$ denote the lattice generated by those lines on the Fermat surface $S$ of degree $m$ that are contained in the set $\mathcal{B}$ of Proposition \ref{Prop:basis}, associated to a root of unity in $\bar\Q$
 that reduces to $\gamma$. In section \ref{s:lines-basis} we have verified for $m \leq 81$ that $M$ has discriminant $\disc (M) = m^{3(m-3)^2}$. For $m > 81$ the discriminant of $M$ is a divisor of a power of $m$ by Corollary \ref{Cor:disc-B}.  
 
 \item
 For each $i$ with $1 \leq i \leq s$ we have 
$\alpha_i \in \F_q$, while $\alpha_i^2+1 = \beta_i^2$ and $\beta_i^q=-\beta_i$. 
In characteristic $p=2$ this implies $\beta_i = \alpha_i+1$, while in odd characteristic it means that $-\beta_i$ is the quadratic conjugate of $\beta_i$ over $\F_q$.  In all cases we have $\beta_i^{q-1} = -1$.
As in section \ref{ss:m},  the line $\li(\alpha_i,\beta_i)$ given by $y = \alpha_i x + \beta_i w$ and $z = \beta_i x + \alpha_i w$
is contained in $(S_{rm})_p$. 
\item
Let $\phi\colon \mathcal S_{rm} \to \mathcal S_m$
be the morphism given by $[x:y:z:w] \mapsto [x^r:y^r:z^r:w^r]$ and set $D_i = \phi (\li(\alpha_i,\beta_i)) \subset (S_m)_p$. Let $L \subset \NS(S_p)$ denote the lattice generated by the image of $M$ and the elements $\sigma \big(D_i\big)$ for all $\sigma \in \mu_m^4/\mu_m$ and all $i$ with $1\leq i \leq s$. 
\end{itemize}

\begin{Result}
\label{Result}
In the above set-up, we have verified with the help of a machine that 
the induced map $M_{\F_\ell} \to L_{\F_\ell}^*$ is injective.
We will comment in \ref{Rem:impl} on some aspects of the implementation.
\end{Result}

 Note that there are two independent reductions involved: the lattice $L$ is contained in the N\'eron-Severi group $\NS(S_p)$ of the reduction of $S$ modulo $p$, while $L_{\F_\ell}$ is the base change of the lattice $L$ over $\Z$ to $\F_\ell$ for a divisor $\ell$ of $m$.
 
\begin{proof}[Proof of Theorem \ref{thm}]
Let $m$ be an integer with $0<m < 100$. If $m>4$ and $(m,6) \neq 1$, then $\NS(S)\otimes \Q$ is not generated by lines by Theorem \ref{Thm:AS}, so certainly $\NS(S)$ is not either. As seen before, for $m \leq 3$ the statement is classical, while for $m=4$ we refer to section \ref{s:4}. We now assume $4<m<100$ and $(m,6)=1$. For $5 \leq m \leq 13$ we could refer to sections \ref{s:5}, \ref{s:7}, \ref{s:11}, and \ref{s:13}, but in any case we can refer to the big table. 
By Corollary \ref{Cor:disc-B},  the discriminant of the lattice $M$, generated by the lines in $\mathcal{B}$ of Proposition \ref{Prop:basis}, is divisible only by primes dividing $m$.

Now we fix the set-up of the table corresponding to the degree $m$.
This involves the supersingular prime $p$ for $\mathcal S$.
Corollary \ref{usefulcor}, applied to $M' = \NS(S)$ and $L' = \NS(S_p)$, shows that $M$ equals $\NS(S)$ thanks to Result \ref{Result}.
We conclude that the N\'eron-Severi lattice $\NS(S)$ is generated by the well-known $3m^2$ lines on $S$ over the $m$-th cyclotomic field, and in fact by those in $\mathcal{B}$. 
\end{proof}

\subsection{Remarks on the implementation}
\label{Rem:impl}

The polynomial $f $ in the table was randomly chosen among the factors of the $m$-th cyclotomic polynomial over $\F_p$.  The pairs $(\alpha_i,\beta_i)$ were chosen randomly among all pairs $(\alpha,\beta) \in \F_q \times \F_{q^2}$ satisfying $\alpha^2+1 = \beta^2$ and $\beta^{q-1} = -1$.  First one pair would be chosen, giving a divisor $D_1$. It was then checked whether the induced map $M_{\F_\ell} \to (L_1)_{\F_\ell}^*$ is injective, where $L_1$ is generated by the image of $M$ in $\NS(S_p)$ and the elements in the orbit of $D_1$ under $\mu_m^4/\mu_m$. In order to save memory, this was not done by writing the entire matrix of intersection numbers between elements of $\mathcal{B}$ on one hand and {\em all} elements of $\mathcal{B}$ and those in the orbit of $D_1$ on the other hand, as there are as many as $m^3$ elements in this orbit. Instead, the kernel of the map $M_{\F_\ell} \to (L_1)_{\F_\ell}^*$ was computed by intersecting the kernel of the map $M_{\F_\ell} \to M_{\F_\ell}^*$ with those of the maps $M_{\F_\ell} \to \langle  \sigma(D_1) \, : \, \sigma \in C \rangle _{\F_\ell}^*$, where $C$ runs through some subsets of $\mu_m^4/\mu_m$ until either the intersection of all kernels was trivial or the union of all subsets $C$ was $\mu_m^4/\mu_4$.  
In order to avoid accidental dependencies, the elements of $C$ were chosen randomly.

 If the computed kernel was not trivial, then a second pair $(\alpha_2,\beta_2)$ was chosen, yielding a divisor $D_2$. The lattice $L_1$ would then be augmented to $L_2$ by also including $D_2$ and the elements in its orbit. The kernel of the new map $M_{\F_\ell} \to (L_2)_{\F_\ell}^*$  would be computed by intersecting the previously computed kernel of $M_{\F_\ell} \to (L_1)_{\F_\ell}^*$ with the kernels of maps 
$M_{\F_\ell} \to \langle  \sigma(D_2) \, : \, \sigma \in C \rangle _{\F_\ell}^*$ for some subsets $C$ of $\mu_m^4/\mu_m$. In all cases this was enough to find a lattice $L$ (namely $L=L_1$ or $L=L_2$) for which $M_{\F_\ell} \to L_{\F_\ell}^*$ is injective.


\end{document}